\documentclass[12pt]{amsart}
\usepackage{amssymb}
\usepackage{amsmath}
\usepackage{longtable}
\newcommand\g{{\mathfrak g}}
\renewcommand\b{{\mathfrak b}}
\newcommand\h{{\mathfrak h}}
\newcommand{\f}{\mathfrak{f}}

\newcommand\gl{\mathfrak{gl}}
\renewcommand{\a}{\mathfrak{a}}

\newcommand\m{\mathfrak m}

\renewcommand\l{\mathfrak l}
\newcommand\n{\mathfrak n}
\newcommand\z{\mathfrak z}
\newcommand\s{\mathfrak s}
\renewcommand{\t}{\mathfrak{t}}

\newcommand\td{\operatorname{tr.deg}}
\newcommand\tr{\operatorname{tr}}

\newcommand\C{\mathbb C}
\newcommand\X{\mathfrak X}

\renewcommand\sl{\mathfrak{sl}}
\newcommand\so{\mathfrak{so}}
\renewcommand\sp{\mathfrak{sp}}
\newcommand\spin{\mathfrak{spin}}
\newcommand\GL{\mathop{\rm GL}\nolimits}

\newcommand\SL{\mathop{\rm SL}\nolimits}

\newcommand\Spin{\mathop{\rm Spin}\nolimits}
\newcommand{\ad}{\mathop{\rm ad}\nolimits}
\newcommand{\Ad}{\mathop{\rm Ad}\nolimits}

\newcommand{\rank}{\mathop{\rm rk}\nolimits}

\newcommand\Int{\mathop{\rm Int}\nolimits}
\newcommand\Aut{\mathop{\rm Aut}\nolimits}
\renewcommand{\Ad}{\mathop{\rm Ad}\nolimits}

\newcommand\quo{/\!/}
\newtheorem{Thm}{Theorem}[section]
\newtheorem{Prop}[Thm]{Proposition}
\newtheorem{Cor}[Thm]{Corollary}
\newtheorem{Lem}[Thm]{Lemma}
\theoremstyle{definition}

\newtheorem{defi}[Thm]{Definition}
\newtheorem{Rem}[Thm]{Remark}
\unitlength=1mm \numberwithin{equation}{section} \oddsidemargin=0cm
\author{Ivan V. Losev}
\title{Computation of  the Cartan spaces of  affine homogeneous spaces}
\date{June 5, 2006}
\thanks{{\it Key words and phrases}: reductive groups, affine homogeneous space, the weight lattice, the Cartan
space} \thanks{{\it 2000 Mathematics Subject Classification.} 14M17,
14R20}
\thanks{The author is supported by RFBR grant  05-01-00988}
\begin{document}
\begin{abstract}
Let $G$ be a reductive algebraic group and $H$ its reductive
subgroup. Fix a Borel subgroup $B\subset G$ and a maximal torus
$T\subset B$. The Cartan space $\a_{G,G/H}$ is, by definition, the
subspace of $\t$, where $\t$ is the Lie algebra of $T$, generated by
the weights of $B$-semiinvariant rational functions on $G/H$.  We
compute the spaces $\a_{G,G/H}$.
\end{abstract}
\maketitle \tableofcontents
\section{Introduction}
The base field $K$ is assumed to be algebraically closed and of
characteristic 0.

Let $G$ be a connected  reductive algebraic group, $T\subset B$ a
maximal torus and a Borel subgroup of $G$, respectively. The
character lattices $\X(T)$ and $\X(B)$ are naturally isomorphic, and
are embedded into $\t^*$. Fix an invariant non-degenerate symmetric
form $(\cdot,\cdot)$ on $\g$ such that $(\xi,\eta)=\tr_V(\xi\eta)$
for some faithful $G$-module $V$ and identify $\g$ with $\g^*$ and
$\t$ with $\t^*$. Note that for any reductive subalgebra $\h\subset
\g$ the restriction of $(\cdot,\cdot)$ to $\h$ is non-degenerate.

Let $X$ be an irreducible $G$-variety. Weights of $B$-semiinvariant
rational functions on $X$ form a sublattice in $\X(T)$. This
sublattice is called the {\it weight lattice} of $X$ and is denoted
by $\X_{G,X}$. The subspace $\a_{G,X}\subset\t^*$ spanned by the
weight lattice is called the {\it Cartan space} of $X$.  The
dimension of $\a_{G,X}$ is called the {\it rank} of $X$ and is
denoted by $\rank_G(X)$.

To justify the terminology consider the case of a symmetric space
$X$. Here $X=G/H$, where $(G^\sigma)^\circ\subset H\subset
G^{\sigma}$,  $\sigma$ is an involution of $G$. Put $\m=\{\xi\in\g|
\sigma(\xi)=-\xi\}$. Choose a maximal commutative subalgebra
$\a\subset \m$ consisting of semisimple elements. It is known that
all such subalgebras are $\Int(\h)$-conjugate. The subspace
$\a\subset\m$ is called the {\it Cartan space} of $G/H$. Choose a
system $\Delta_+$ of positive roots in $\a$.  Let $\b$ be a Borel
subalgebra of $\z_\g(\a)\oplus\bigoplus_{\alpha\in
\Delta_+}\g^{\alpha}$ and $\t$ a Cartan subalgebra of
$\b\cap\z_\g(\a)$.  It turns out (see, for example, \cite{Vust2})
that $\a=\a_{G,G/H}$.

The goal of this article is to compute the subspaces
$\a_{G,G/H}\subset\t$, where $H$ is a reductive subgroup of $G$.
Note that $\a_{G,G/H}\subset \a_{G,G/H_1}$ for any normal subgroup
$H_1\subset H$ because there is a dominant $G$-equivariant morphism
$G/H_1\rightarrow G/H$. Moreover, since $K(G/H^\circ)$ is an
algebraic extension of $K(G/H)$, $\a_{G,G/H}=\a_{G,G/H^\circ}$. Thus
the subspace  $\a_{G,G/H}$ depends only on the pair $(\g,\h)$  and
we put $\a(\g,\h)=\a_{G,G/H}$. The main idea of the computation is
to notice that the space $\a(\g,\h)$ depends only on a certain
"essential" ideal of $\h$. There is a precise definition:

\begin{defi}\label{Def:1.1}
A reductive subalgebra  $\h\subset \g$ is called {\it essential}, if
for any proper ideal $\h_1\subset \h$ the inclusion
$\a(\g,\h)\subset \a(\g,\h_1)$ is strict.
\end{defi}

To state our result we also we need a quite standard notion of an
indecomposable subalgebra.

\begin{defi}\label{Def:1.2}
Suppose that a subalgebra $\h\subset\g$ is such that there exist
ideals $\g_1,\g_2\subset\g$ and reductive subalgebras $\h_1\subset
\g_1,\h_2\subset \g_2$ with $\g=\g_1\oplus\g_2$,
$\h=\h_1\oplus\h_2$. In this situation we say that $(\g,\h)$ is the
direct sum of the pairs $(\g_1,\h_1)$, $(\g_2,\h_2)$ and write
$(\g,\h)=(\g_1,\h_1)\oplus (\g_2,\h_2)$. If the pair $(\g,\h)$ can
not be decomposed into a non-trivial direct sum, then $\h$ is said
to be  an {\it indecomposable} subalgebra.
\end{defi}

Clearly, $\a(\g_1\oplus\g_2,\h_1\oplus\h_2)=\a(\g_1,\h_1)\oplus
\a(\g_2,\h_2)$.

Now we are ready to state our main result.

\begin{Thm}\label{Thm:1.3}
Let $G,H,\g,\h$ be such as above. Then
\begin{itemize}
\item[(a)] There exists a unique ideal $\h^{ess}\subset \h$ such
that $\h^{ess}$ is an essential subalgebra in $\g$ and
$\a(\g,\h)=\a(\g,\h^{ess})$. This ideal is called the {\it essential
part} of $\h$. The ideal $\h^{ess}$ is maximal (with respect to an
inclusion) among the ideals of $\h$, that are essential subalgebras
of $\g$. If $(\g,\h)=(\g_1,\h_1)\oplus (\g_2,\h_2)$, then
$\h^{ess}=\h_1^{ess}\oplus \h_2^{ess}$.
\item[(b)] All non-trivial essential semisimple indecomposable subalgebras  up to  conjugacy in $\Aut(\g)$
are precisely those given in Table~\ref{Tbl:1.4}.
\item[(c)] In this part we classify indecomposable essentail
nonsemisimple subalgebras and show how to compute the Cartan spaces
for them.
\begin{enumerate}
\item
Let $\h$ be an indecomposable essential not semisimple subalgebra of
$\g$. Then the pair $([\g,\g],[\h,\h])$ is the direct sum  of some
copies of
 pairs 1,2 (the latter for $k\neq n/2$),10,19 from
Table~\ref{Tbl:1.4},
 and
the pairs   $(\sl_{2n+1},\sp_{2n}), (\sl_{2n+1},\sl_{n+1})$. In
parts c2-c5 we assume that $[\h,\h]$ has the indicated form.
\item  $\z(\z_\g(\h))= \z(\z_\g([\h,\h]))$.
\item Let $\z_1\subset \z(\z_\g([\h,\h]))$ be an algebraic
subalgebra. If $\z_1\oplus [\h,\h]$ is indecomposable, then it is
essential.
 \item The spaces $\z:=\z(\z_\g([\h,\h]))$ and
$\a(\g,[\h,\h])/\a(\g,[\h,\h]\oplus\z)$ are dual to each other. The
duality is given by the rule formulated in Remark~\ref{Rem:1.5}.
\item The subspace $\a(\g,\h)\subset \a(\g,[\h,\h])$ is the inverse
image under the projection $\a(\g,[\h,\h])\rightarrow
\a(\g,[\h,\h])/\a(\g,[\h,\h]\oplus\z)$ of the annihilator of the
subspace  $\z(\h)\subset \z$.
\end{enumerate}
\item[(d)] Let $\sigma$ be an automorphism of $\g$. Then the
relation between $\a(\g,\h)$ and $\a(\g,\sigma(\h))$ is given by
Remark~\ref{Rem:1.7}.
\end{itemize}
\end{Thm}
\setlongtables
\begin{longtable}{|c|c|c|c|}
\caption{Semisimple indecomposable essential subalgebras
$\h\subset\g$}\label{Tbl:1.4}\\\hline
N&$\g$&$\h$&$\a(\g,\h)$\\\endfirsthead\hline
N&$\g$&$\h$&$\a(\g,\h)$\\\endhead\hline 1&$\sl_n,n\geqslant
2$&$\sl_k, \frac{n+2}{2}\leqslant k\leqslant
n$&$\langle\pi_i,\pi_{n-i}; i\leqslant n-k\rangle$\\\hline
2&$\sl_n,n\geqslant 4$&$\sl_k\times\sl_{n-k}, \frac{n}{2}\leqslant
k\leqslant n-2$&$\langle\pi_i+\pi_{n-i},\pi_k,\pi_{n-k}; i<
n-k\rangle$
\\\hline 3&$\sl_{2n},n\geqslant
2$&$\sp_{2n}$&$\langle\pi_{2i}; i\leqslant n-1\rangle$\\\hline
4&$\sp_{2n},n\geqslant 2$&$\sp_{2k},\frac{n+1}{2}\leqslant
k\leqslant n$&$\langle\pi_i; i\leqslant 2(n-k)\rangle$\\\hline
5&$\sp_{2n},n\geqslant 2$&$\sp_{2k}\times\sp_{2(n-k)},
\frac{n}{2}\leqslant k< n $&$\langle\pi_{2i}; i\leqslant
n-k\rangle$\\\hline 6&$\sp_{2n},n\geqslant 4
$&$\sp_{2n-4}\times\sl_2\times\sl_2$&$\langle
\pi_2,\pi_4,\pi_1+\pi_3\rangle$\\\hline
7&$\sp_6$&$\sl_2\times\sl_2\times\sl_2$&$\langle\pi_2,\pi_1+\pi_3\rangle$\\\hline
 8&$\so_n,n\geqslant
7$&$\so_k, \frac{n+2}{2}\leqslant k\leqslant
n$&$\langle\pi_i,i\leqslant n-k\rangle$\\\hline
9&$\so_{4n},n\geqslant 2$&$\sl_{2n}$&$\langle\pi_{2i}; i\leqslant
n\rangle$\\\hline 10&$\so_{4n+2},n\geqslant 2
$&$\sl_{2n+1}$&$\langle \pi_{2i},\pi_{2n+1}; i\leqslant
n\rangle$\\\hline
11&$\so_9$&$\spin_7$&$\langle\pi_1,\pi_4\rangle$\\\hline
12&$\so_{10}$&$\spin_7$&$\langle
\pi_1,\pi_2,\pi_4,\pi_5\rangle$\\\hline
13&$\so_7$&$G_2$&$\langle\pi_3\rangle$\\\hline
14&$\so_8$&$G_2$&$\langle \pi_1,\pi_3,\pi_4\rangle$\\\hline
15&$G_2$&$A_2$&$\langle\pi_1\rangle$\\\hline
16&$F_4$&$B_4$&$\langle\pi_1\rangle$\\\hline
17&$F_4$&$D_4$&$\langle\pi_1,\pi_2\rangle$\\\hline
18&$E_6$&$F_4$&$\langle\pi_1,\pi_5\rangle$\\\hline
19&$E_6$&$D_5$&$\langle\pi_1,\pi_5,\pi_6\rangle$\\\hline
20&$E_6$&$B_4$&$\langle\pi_1,\pi_2,\pi_4,\pi_5,\pi_6\rangle$\\\hline
21&$E_6$&$A_5$&$\langle \pi_1+\pi_5,\pi_2+\pi_4,
\pi_3,\pi_6\rangle$\\\hline
22&$E_7$&$E_6$&$\langle\pi_1,\pi_2,\pi_6\rangle$\\\hline
23&$E_7$&$D_6$&$\langle \pi_2,\pi_4,\pi_5,\pi_6\rangle$\\\hline
24&$E_8$&$E_7$&$\langle \pi_1,\pi_2,\pi_3,\pi_7\rangle$\\\hline
25&$\h\times\h$&$\h$&$\langle\pi^*_i+\pi'_i;i\leqslant
\rank\h\rangle$\\\hline 26&$\sp_{2n}\times\sp_{2m},
m>n>1$&$\sp_{2n-2}\times\sl_2\times\sp_{2m-2}$&$\langle\pi_2,\pi_2',\pi_1+\pi_1'\rangle$\\\hline
27& $\sp_{2n}\times \sl_2,n>1$&$\sp_{2n-2}\times \sl_2$&$\langle
\pi_2, \pi_1+\pi_1'\rangle$\\\hline
\end{longtable}

Most of the notation used in the table is explained in
Section~\ref{SECTION_Not_Conv}. In rows 25-27  $\pi_i$ (resp.,
$\pi_i'$) denotes the fundamental weight of the first (resp.,
second) ideal.

\begin{Rem}\label{Rem:1.5}
This remark explains how to compute the spaces $\a(\g,\h)$ in the
case, when $\h$ is an essential not semisimple subalgebra.

Let $\g=\z(\g)\oplus\g_1\oplus\ldots\oplus\g_k$ be the decomposition
into  the direct sum of the center and simple (non-commutative)
ideals. Recall that in the interesting for us cases the pair
$([\g,\g],[\h,\h])$ is the direct sum of  pairs NN 1,2 ($k\neq
n/2$),10,19 from Table~\ref{Tbl:1.4} and the pairs
$(\sl_{2n+1},\sp_{2n}),(\sl_{2n+1},\sl_{n+1})$. Put
$\h_i=[\h,\h]\cap\g_i$, $\z_0=\z(\g), \z_i=\z(\z_{\g_i}(\h_i)),
i=\overline{1,k}$.  It is clear that
$\z=\z(\z_\g([\h,\h]))=\z_0\oplus\z_1\oplus\ldots\oplus\z_k$ and
$$\a(\g,[\h,\h])/\a(\g,[\h,\h]\oplus \z
)=\a(\z(\g),0)\oplus\bigoplus_{i=1}^k
\a(\g_i,\h_i)/\a(\g_i,\h_i\oplus \z_i).$$

 The identification of $\z$ with
$(\a(\g/[\h,\h])/\a(\g,[\h,\h]\oplus \z))^*$ is established as
follows. Let $G$ be a simply connected group corresponding to $\g$
and $H$ a connected subgroup of $G$ corresponding to $\h$. An
element $x\in\z$ is identified with a unique linear function
$\alpha_x\in\a(\g,[\h,\h])^*$ such that the equality
$xv=\alpha_x(\lambda)v$ holds for any highest weight $\lambda$  and
any  $v\in V(\lambda)^{(H,H)}$. Using the Frobenius reciprocity, we
see the conditions
\begin{enumerate}
\item $V(\lambda)^{(H,H)}\neq 0$,
\item $V(\lambda^*)^{(H,H)}\neq 0$,
\item $\lambda$ is a highest weight of $K[G/(H,H)]$
\end{enumerate}
are equivalent.

We show in Section~\ref{SECTION_computation} that for all six pairs
$([\g,\g],[\h,\h])$ listed above  $\alpha_x$ is well-defined and
annihilates $\a(\g,[\h,\h]\oplus\z)$.

Now we describe the correspondence $x\mapsto \alpha_x$ explicitly.

It follows from the definition that for  $x\in\z_i,
i=\overline{0,k},$ the function $\alpha_x$ lies in $\a(\z(\g),0)^*$
for $i=0$ and in $(\a(\g_i,\h_i)/\a(\g_i,\h_i\oplus \z_i))^*$ for
$i>0$. For $i=0$ the map $x\in\z_0\mapsto\alpha_x$ is the natural
identification $\z_0\cong\a(\z_0,0)^*$. For $i>0$ $\alpha_x$
coincides with the analogous function defined for the pair
$(\g_i,\h_i)$. Now its enough to find all $\a(\g_i,\h_i\oplus\z_i)$
and to identify $(\a(\g_i,\h_i)/\a(\g_i,\h_i\oplus \z_i))^*$ with
$\z_i$. All spaces $\z_i$ are one-dimensional, so  it is enough to
find a highest weight $\lambda$ of $K[G_i/H_i]$ not lying in
$\a(\g_i,\h_i\oplus \z_i)$ (here $G_i,H_i$ denote the connected
reductive groups corresponding  to $\g_i$, $\h_i$) and the number
$\alpha_x(\lambda_i)$. This information is given in Table
~\ref{Tbl:1.6}. We denote by $\pi_i^\vee$ the dual fundamental
weight.
\end{Rem}

\begin{longtable}{|c|c|c|c|c|}
\caption{The identification of $\z_i$ with
$(\a(\g_i,\h_i)/\a(\g_i,\h_i\oplus\z_i))^*$}\label{Tbl:1.6}\\\hline
$(\g_i,\h_i)$&$x\in\z_i$&$\lambda$&$\alpha_x(\lambda)$&$\a(\g_i/\h_i\oplus
\z_i)$\\\hline $(\sl_{n},\sl_k),
k>\frac{n}{2}$&$\pi_{n-k}^\vee$&$\pi_1$&$\frac{k}{n}$&$\{\sum_{i=1}^{n-k}(x_i\pi_i+x_{n-i}\pi_{n-i});$\\
&&&&$\sum_{i=1}^{n-k}i(x_i-x_{n-i})=0\}$\\\hline
$(\sl_n,\sl_{k}\times\sl_{n-k}),k>
\frac{n}{2}$&$\pi_{n-k}^{\vee}$&$\pi_{n-k}$&$\frac{(n-k)k}{n}$&$\langle
\pi_i+\pi_{n-i}; i\leqslant n-k\rangle$\\\hline
$(\sl_{2n+1},\sp_{2n})$&$\pi_1^\vee$&$\pi_1$&$\frac{2n}{2n+1}$&$\{\sum_{i=1}^{2n}x_i\pi_i;$\\
&&&&$\sum_{i=0}^{n-1}(n-i)x_{2i+1}-\sum_{i=1}^nix_{2i}=0\}$\\\hline
$(\so_{4n+2},\sl_{2n+1})$&$\pi_{2n+1}^\vee$&$\pi_{2n+1}$&$n+\frac{1}{2}$&$\langle
\pi_{2i},\pi_{2n}+\pi_{2n+1}; i\leqslant n-1\rangle$\\\hline
$(E_6,D_5)$&$\pi_1^\vee$&$\pi_1$&$\frac{4}{3}$&$\langle
\pi_1+\pi_5,\pi_6\rangle$\\\hline
\end{longtable}

\begin{Rem}\label{Rem:1.7}
This remark explains the relation between the spaces $\a(\g,\h)$ and
$\a(\g,\sigma(\h))$ for $\sigma\in \Aut(\g)$. First of all, we note
that there exist $g\in \Int(\g)$ such that $\b,\t$ are
$\sigma_1$-stable, where $\sigma_1=g\sigma$. Changing $\sigma$ by
$\sigma_1$, we may assume that $\sigma(\b)=\b,\sigma(\t)=\t$. In
this case $\a(\g,\sigma(\h))=\sigma(\a(\g,\h))$.

This can be seen, for example, in the following way. By the
Frobenius reciprocity, $\a(\g,\h)$ (respectively,
$\a(\g,\sigma(\h))$) is generated by all heighest weights $\lambda$
such that $V(\lambda)^{*\h}\neq 0$ (resp.,
$V(\lambda)^{*\sigma(\h)}\neq 0$). But $\dim
V(\lambda)^{*\sigma(\h)}=\dim V(\sigma(\lambda))^{*\h}$.

Note that for all subalgebras from
Tables~\ref{Tbl:1.4},\ref{Tbl:1.6} except NN8 ($n=8,k=7$), 9, 25 any
subalgebra $\Aut(\g)$-conjugate to $\h$  is also
$\Int(\g)$-conjugate. In case $8$ (resp., 9,25) the class of
$\Aut(\g)$-conjugacy of $\h$ is the union of  3 (resp. 2,$\#
\Aut(\h)/\Int(\h)$) classes of $\Int(\g)$-conjugacy.
\end{Rem}

\section{Basic properties of essential
subalgebras}\label{SECTION_essential} In this subsection $G$ is a
connected reductive group, $H$ its reductive subgroup.

For an irreducible  $G$-variety $X$ we put
\begin{align}\label{eq:2.1}
&L_{G,X}=Z_G(\a_{G,X}).\\\label{eq:2.2} &L_{0\,G,X}=\cap_{\chi\in
\X_{G,X}}\ker\chi.
\end{align}

In (\ref{eq:2.2}) we consider $\X_{G,X}$ as a sublattice in
$\X(L_{G,X})$.

We need the following crucial fact proved (in greater generality) by
F.Knop in~\cite{Knop1}, Satz 8.1, Korollar 8.2.

\begin{Prop}\label{Prop:2.1}
Let $X$ be a quasiaffine smooth $G$-variety. Then $L_{0\,G,X}$ is
the stabilizer in general position for the action
 $G:T^*X$.
\end{Prop}

We recall that, by definition, a subgroup  $G_0\subset G$, where $G$
is an arbitrary algebraic group, is called the {\it stabilizer in
general position} (s.g.p.) for an irreducible $G$-variety $Y$ if for
some open subset $Y^0\subset Y$ the stabilizer of any $y\in Y^0$ is
conjugate to $G_0$ in $G$. The Lie algebra of the s.g.p. is called
the {\it stable subalgebra in general position} (s.s.g.p.) for the
action $G:Y$. If $Y$ is a vector space, $G$ is reducitive and the
action $G:Y$ is linear, then there exists the s.g.p for this action,
see~\cite{VP}, $\S7$.

\begin{Rem}\label{Rem:2.2}
The space $\a_{G,X}$ is not necessarily determined uniquely by the
s.s.g.p for the action $G:T^*X$. For example, let $G=\SL_n,
X_1=K^n\oplus K^n,X_2=K^n\oplus K^{n*}, X_3=K^{n*}\oplus K^{n*}$,
where $K^n$ is the space of  the tautological representation of
$\SL_n$. Then $\a_{G,X_1}=\langle \pi_1,\pi_2\rangle,
\a_{G,X_2}=\langle \pi_1,\pi_{n-1} \rangle, \a_{G,X_3}=\langle
\pi_{n-1},\pi_{n-2}\rangle$. But $T^*X_1\cong T^*X_2\cong T^*X_3$.
\end{Rem}

When $X=G/H$, one obtains the following corollary from
Proposition~\ref{Prop:2.1}.
\begin{Cor}\label{Cor:2.3}
The subalgebra $\l_{0\,G,G/H}\subset \g$ is conjugate to the s.s.g.p
for the $\h$-module $\g/\h\cong (\g/\h)^*\cong \h^\perp$.
\end{Cor}

In the sequel we need two simple lemmas.

\begin{Lem}\label{Lem:2.4}
Let $\h$ be a reductive Lie algebra, $V$ an $\h$-module (see the
conventions after the index of notation), $\l_0$ the corresponding
s.s.g.p and $\h_1$ a reductive subalgebra in $\h$. Then the
subalgebra $\Ad(h)\l_0\cap \h_1$ for $h\in \Int(\h)$ in general
position is the s.s.g.p. for the $\h_1$-module $V$. In particular,
if $\h_1\subset \h$ is an ideal, then $\l_0\cap \h_1$ is the s.s.g.p
for the $\h_1$-module $V$.
\end{Lem}
\begin{proof}
This follows directly from the definition of the s.s.g.p.
\end{proof}

\begin{Lem}\label{Lem:2.5}
Let $\g$ be a reductive Lie algebra, $\h$ a reductive subalgebra in
$\g$ and $\h_1$ an ideal in $\h$. Then the $\h_1$-module $\g/\h_1$
is the direct sum of a trivial $\h_1$-module and  the $\h_1$-module
$\g/\h$.
\end{Lem}
\begin{proof}
This follows from the fact that $\h/\h_1$ is a trivial
$\h_1$-module.
\end{proof}

The main result of this section is the following
\begin{Prop}\label{Prop:2.6}
 The ideal $\h^{ess}\subset \h$, defined in part (a) of
Theorem~\ref{Thm:1.3} exists and is uniquely determined. Suppose
$\l_0$ is the s.s.g.p for the $\h$-module $\g/\h$. Then $\h^{ess}$
coincides with the ideal in $\h$ generated by $\l_0$.
\end{Prop}
\begin{proof}
Denote by $\h_1$ the ideal in $\h$ generated by $\l_0$. By
Lemma~\ref{Lem:2.4}, $\l_0$ is the s.s.g.p for the $\h_1$-module
$\g/\h$, and hence, by Lemma~\ref{Lem:2.5}, also for  $\g/\h_1$.
Since $\a(\g,\h)\subset \a(\g,\h_1)$  and
$\dim\a(\g,\h)=\rank\g-\rank\l_0=\dim\a(\g,\h_1)$, we obtain
$\a(\g,\h)=\a(\g,\h_1)$.

Next we show that any ideal  $\h_2\subset \h$ with
$\a(\g,\h_2)=\a(\g,\h)$ contains $\h_1$. Let $H_2$ be a connected
subgroup of $H$ corresponding to $\h_2$. Then
$\l_{0\,G,G/H_2}=\l_{0\,G,G/H}$. Corollary~\ref{Cor:2.3} and
Lemma~\ref{Lem:2.5} imply that $\l_{0}$ and the s.s.g.p for the
$\h_2$-module $\g/\h_2$ are $G$-conjugate. It follows from
Lemma~\ref{Lem:2.4} that $\l_0\subset \h_2$. By the definition of
$\h_1$, $\h_1\subset \h_2$.

If $\h_2$ is an ideal of $\h_1$ satisfying
$\a(\g,\h_2)=\a(\g,\h_1)$, then $\h_1\subset \h_2$ because $\h_2$ is
an ideal of $\h$. Thus $\h_1$ is an essential subalgebra of $\g$.

The equality $\a(\g,\h_1)=\a(\g,\h)$ implies that if $\h_2$ is an
ideal of $\h$ containing $\h_1$, then $\h_2$ is essential iff
$\h_2=\h_1$. This completes the proof.
\end{proof}

\begin{Lem}\label{Lem:2.9}
Let $\h$ be a reductive Lie algebra,  $V_1$ an $\h$-module, $V_2$ an
$\h$-submodule in $V_1$, $\l_i$ the s.s.g.p. for $V_i, i=1,2$. Then
$\l_1$ is $\Int(\h)$-conjugate to the s.s.g.p. for the $\l_2$-module
$V_1/V_2$. In particular, if the ideal of $\h$ generated by $\l_1$
coincides with $\h$, then the same property holds for $\l_2$.
\end{Lem}
\begin{proof}
Choose an invariant complement $V_2'$ to $V_2$ in $V_1$. Then, up to
the $\Int(\h)$-conjugacy, $\l_1=\h_{v_1+v_2}=(\h_{v_2})_{v_1},
\l_2=\h_{v_2}$ for $v_1\in V_2',v_2\in V_2$ in general position.
\end{proof}

\begin{Cor}\label{Cor:2.7}
If $\h_1\subset\h_2$ are reductive subalgebras of $\g$, then
$\h_1^{ess}\subset \h_2^{ess}$.
\end{Cor}
\begin{proof}
Let $\l_{0i}$ be the s.s.g.p. for the $\h_i$-module $\g/\h_i,
i=1,2$. It follows directly from the definition of s.s.g.p's that
the s.s.g.p for the $\h_1$-module $\g/\h_2$ is
$\Int(\h_2)$-conjugate to a subalgebra of $\l_{02}$. Using
Lemma~\ref{Lem:2.9}, we see that $\l_{01}$ is $\Int(\h_2)$-conugate
to a subalgebra of $\l_{02}$. It remains to apply
Proposition~\ref{Prop:2.6}.
\end{proof}

\begin{Cor}\label{Cor:2.12}
Let  $\h$ be an essential subalgebra in $\g$. Then the following
assertions hold
\begin{enumerate}
\item Let $\g_0$ be a reductive subalgebra of $\g$, $\h_0$ a reductive subalgebra of
$\g_0$ and $\varphi:\h\rightarrow \h_0$ an epimorphism of Lie
algebras such that the $\h$-module $\g_0/\h_0$ is a submodule in
$\g/\h$. Then $\h_0$ is an essential subalgebra in $\g_0$.
\item The projection of  $\h$ to any ideal  $\g_1\subset\g$
is an essential subalgebra in $\g_1$.
\item Let $\h_0$ be an ideal in $\h$. Then $\h/\h_0$ is an essential subalgebra in $\n_\g(\h_0)/\h_0$.
\end{enumerate}
\end{Cor}
\begin{proof}
The first assertion is deduced from Lemma~\ref{Lem:2.9} and the
characterization of essential subalgebras given in
Proposition~\ref{Prop:2.6}. The second and the third assertions are
special cases of the first one.
\end{proof}

\begin{Cor}\label{Cor:2.13}
Suppose $\h$ does not contain ideals of $\g$. Let $\a$ be a simple
non-commutative ideal in $\h$ such that there exist different simple
ideals $\g_1,\ldots,\g_k\subset \g$ such that the projection of $\a$
to $\g_i$ is nonzero for all $i$. Then
\begin{enumerate}
\item  The $\h$-module
$\g/\h$ contains an $\h$-submodule  isomorphic to $\a^{\oplus k-1}$.
\item
Let $\h\subset\g$ be an essential subalgebra. Then $k\leqslant 2$.
\end{enumerate}
\end{Cor}
\begin{proof}
The multiplicity of the $\h$-module $\a$ in the $\h$-modules $\g_i,
i=\overline{1,k},$ and $\h$  equals 1. Thus the multiplicity of $\a$
in $\g$ is not less than $k$. This proves the first assertion.

Assume now that $\h$ is essential but $k\geqslant 3$. Then the
$\h$-module $\g/\h$ contains a submodule $\a\oplus \a$. The s.s.g.p.
of the $\a$-module $\a\oplus \a$ is trivial. Indeed, the s.s.g.p.
for the $\a$-module $\a$ is a Cartan subalgebra $\t_0\subset \a$.
The s.s.g.p. of the $\t_0$-module $\a$ is trivial because $\t_0$ is
commutative and the representation of $\t_0$ in $\a$ is effective.
We see that the s.s.g.p. for the $\h$-module $\a\oplus \a$ projects
trivially to $\a$. By Lemma~\ref{Lem:2.9}, the same holds for the
s.s.g.p for $\g/\h$. Proposition~\ref{Prop:2.6} implies $\h$ is not
essential.
\end{proof}

\begin{Rem}\label{Rem:2.8}
The ideal generated by a  subalgebra $\f\subset \h$ coincides with
the direct sum of the projection of $\f$ to $\z(\h)$ and all simple
ideals $\h_0\subset \h$ such that the projection of $\f$ to $\h_0$
is nonzero.
\end{Rem}

\begin{proof}[Proof of part (a) of Theorem~\ref{Thm:1.3}]
At first, we prove that $\h^{ess}$ contains any ideal $\h_1\subset
\h$ such that $\h_1$ is an essential subalgebra of $\g$. Let $\l_0$
be the s.s.g.p for the $\h$-module $\g/\h$. The subalgebra $\l_0\cap
\h_1$ is the s.s.g.p for the $\h_1$-module $\g/\h_1$
(Lemmas~\ref{Lem:2.4},\ref{Lem:2.5}). Apply
Proposition~\ref{Prop:2.6} to the pair $(\g,\h_1)$. We see that
there is no proper ideals of $\h_1$ containing  $\l_0\cap \h_1$.
Thus $\h_1$ is an ideal of $\h$ generated by $\l_0\cap \h_1$. Since
the ideal $\h^{ess}\subset \h$ is generated by $\l_0$, $\h_1\subset
\h^{ess}$.

Now let $(\g,\h)=(\g_1,\h_1)\oplus (\g_2,\h_2)$. Denote by
$\l_{0i},i=1,2,$ the s.s.g.p for the $\h_i$-module $\g_i/\h_i$. Then
$\l_{01}\oplus \l_{02}$ is the s.s.g.p for the $\h$-module $\g/\h$.
Proposition~\ref{Prop:2.6} implies $\h^{ess}=\h_1^{ess}\oplus
\h^{ess}_2$.

Proposition~\ref{Prop:2.6} implies all remaining claims.
\end{proof}

To make the description of essential subalgebras more convenient we
need the  notion of the {\it saturation} of an essential subalgebra.

\begin{defi}\label{Def:2.10}
 Denote by $\widetilde{\h}$ the inverse image of $\z(\n_\g(\h)/\h)$
under the natural projection $\n_{\g}(\h)\rightarrow \n_\g(\h)/\h$.
By definition, the {\it saturation} of $\h$ is
$\widetilde{\h}^{ess}$. We denote the saturation of $\h$ by
$\h^{sat}$. We say that
 $\h$ is {\it saturated} if $\h=\h^{sat}$.
\end{defi}

Suppose now $\h$ is essential. Then $\h\subset \h^{sat}$, by
Corollary~\ref{Cor:2.7}. Clearly, $\h^{sat}$ is the direct sum of
$\h$ and some commutative ideal. Let $\l_0,\l_0^1$ be the s.s.g.p's
for the $\h$-module $\g/\h$ and the $\h^{sat}$-module $\g/\h^{sat}$,
respectively. By Lemmas~\ref{Lem:2.4},\ref{Lem:2.5},  $\l^1_0\cap
\h=\l_0$. This allows us to reduce the classification of arbitrary
essential subalgberas to the classification of saturated ones
(together with the corresponding s.s.g.p). Namely, let $\h$ be a
saturated subalgebra. Then the following conditions are equivalent
\begin{enumerate} \item an ideal $\h_1\subset \h$ containing
$[\h,\h]$ is an essential subalgebra of $\g$. \item the subalgebra
$\l_0\cap\h_1$ which is the s.s.g.p for the $\h_1$-module $\g/\h_1$
projects non-trivially to any simple non-commutative ideal of
$\h$.\end{enumerate} Indeed, $\l_0\cap\h_1$ projects surjectively
onto $\z(\h_1)=\h_1\cap\z(\h)$ because $\l_0$ projects surjectively
onto $\z(\h)$. The equivalence of the two conditions follows now
from Remark~\ref{Rem:2.8}. Thus there is a one-to-one correspondence
between the set of all essential subalgebras $\h_1\subset \g$ with
$\h_1^{sat}=\h$ and the set of all algebraic subalgebras
$\t_1\subset \l_0/([\h,\h]\cap\l_0)\cong \z(\h)$ such that the
inverse image of $\t_1$ in $\l_0$ has a non-zero projection to any
simple non-commutative ideal of $[\h,\h]$. This correspondence is
given by $\h_1\mapsto \h_1/[\h,\h]$.

 Note also that if
$(\g,\h)=(\g_1,\h_1)\oplus (\g_2,\h_2)$, then
$\h^{sat}=\h_1^{sat}\oplus \h_2^{sat}$.

We say that a subalgebra $\h$ of $\g$ is {\it initial} if $\h$ is
 saturated and indecomposable. The classification of essential
 subalgebras
reduces to the classification of initial subalgebras and the
calculation of the subalgebras $\l_0$ for all initial subalgebras.
Clearly, any initial subalgebra $\h\subset \g$ contains $\z(\g)$.
Thus it is enough to classify initial subalgebras in semisimple
algebras.

\section{Indeces of subalgebras and modules}\label{SECTION_indeces}
We use notation of  Section~\ref{SECTION_essential}.

In  the previous section we have noticed the close connection
between essential subalgebras and s.s.g.p. for certain modules. In
the study of s.s.g.p's the notion of the {\it index} of a module
over a simple Lie algebra plays a great role. Let $\h$ be a
reductive Lie algebra, $V$ an $\h$-module (see conventions on
modules after the index of notation). Define an invariant symmetric
bilinear  form $(\cdot,\cdot)_V$ on $\h$ by the formula
$(x,y)_V=\tr_V(xy)$. If the module $V$ is effective, then
$(\cdot,\cdot)_V$ is non-degenerate. Suppose that $\h$ is simple.
Then an invariant symmetric bilinear form is determined uniquely up
to the multiplication by a constant (see, for example,~\cite{AEV}).
The {\it index} of the $\h$-module $V$ is, by definition,
$\frac{(x,y)_V}{(x,y)_\h}$. The last fraction does not depend on the
choice of $x,y$ with $(x,y)_\h\neq 0$. We denote the index by
$l_\h(V)$.

The following proposition is a straightforward generalization of
results from \cite{AEV}, \cite{Elash1}.

\begin{Prop}\label{Prop:3.1}
Let $\h$ be a reductive Lie algebra, $V$ an
 effective $\h$-module. Let $\h=\z(\h)\oplus \h_1\ldots\oplus\h_k$
be the decomposition of $\h$ into the direct sum of the center and
simple non-commutative ideals. Then:
\begin{enumerate}
\item If the s.s.g.p for the $\h$-module $V$
is non-trivial, then $l_{\h_i}(V)\leqslant 1$ for some $i$.
\item If $l_{\h_i}(V)\geqslant 1$ for all $i$, then the s.s.g.p is
contained in $\bigoplus_{i, l_{\h_i}(V)=1}\h_i$.
\end{enumerate}
\end{Prop}
\begin{proof}
Let $\l_0$ be the s.s.g.p for the $\h$-module $V$. This is an
algebraic subalgebra of $\h$. If $\l_0$ contains a nilpotent
element,  then the both assertions follow from the proof of the main
theorem in~\cite{AEV}. So we may assume that $l_{\h_i}(V)\geqslant
1$ for all $i$ and that $\l_0$ contains a rational semisimple
element $x$ (the word "rational" means "having rational eigenvalues
on any $\h$-module"). Put $x=x_0+x_1+\ldots +x_k$, where
$x_0\in\z(\h), x_i\in\h_i,$ for $i=\overline{1,k}$. It follows from
Lemma 3 in~\cite{AEV} that
\begin{equation*}\sum_{i=0}^k(x_i,x_i)_V=(x,x)_V\leqslant
(x,x)_\h=\sum_{i=1}^k (x_i,x_i)_{\h_i}.
\end{equation*}
This implies  the both assertions of the proposition  (note that
$(x_0,x_0)_V\neq 0$ for any rational $x_0\in\z(\h)$ because $V$ is
effective).
\end{proof}

Our next objective is to describe all reductive subalgebras
$\h\subset\g$ such that $l_{\h_1}(\g/\h)\leqslant 1$ for some simple
ideal $\h_1\subset\h$. To do this we need the notion of the Dynkin
index.

Until the end of the section $\g,\h$ are simple Lie algebras,
$\iota$ is an embedding $\h\hookrightarrow \g$. We fix an invariant
non-degenerate symmetric bilinear form $K_\g$ on $\g$ such that
$K_\g(\alpha^{\vee},\alpha^{\vee})=2$ for a root
$\alpha\in\Delta(\g)$ of the maximal length. Analogously define a
form $K_\h$ on $\h$. Recall that the  {\it Dynkin index} of the
embedding $\iota:\h\hookrightarrow \g$  is, by definition,
$K_\g(\iota(x),\iota(x))/K_\h(x,x)$ (the last fraction does not
depend on the choice of $x\in\h$ such that $K_\h(x,x)\neq 0$).
Abusing the notation, we denote the Dynkin index of $\iota$ by
$i(\h,\g)$. It turns out that $i(\h,\g)$ is a positive integer
(see~\cite{Dynkin}).

By $k_\g$ we denote the number $(\alpha^\vee,\alpha^\vee)_\g$ for a
root  $\alpha\in\g$ of  maximal length. Analogously define $k_\h$.
It follows immediately from the definition that
\begin{equation*}
k_\g=\sum_{\beta\in\Delta(\g)}\langle
\beta,\alpha^{\vee}\rangle^2=8+\#\{\beta\in\Delta(\g)| \beta\neq
\pm\alpha, (\alpha,\beta)\neq 0\}.
\end{equation*}

The numbers $k_\g$ for all simple Lie algebras are given in
Table~\ref{Tbl:3.2}.

\begin{longtable}{|c|c|c|c|c|c|c|c|c|c|}
\caption{$k_\g$. }\label{Tbl:3.2}\\\hline
$\g$&$A_l$&$B_l$&$C_l$&$D_l$&$E_6$&$E_7$&$E_8$&$F_4$&$G_2$\\\hline
$k_\g$& $4l+4$& $8l-4$&$4l+4$&$8l-8$&48&72& 120&36&16\\\hline
\end{longtable}

Analyzing possible embeddings between simple Lie algebras we obtain
the following
\begin{Lem}\label{Lem:3.3}
 $k_\h<k_\g$.
\end{Lem}

By the definition of the the index $l_{\h}(\g/\h)$, we get
\begin{equation}\label{eq:3.1}
l_\h(\g/\h)=\frac{i(\h,\g)k_\g}{k_\h}-1.\end{equation} It follows
from (\ref{eq:3.1}) and Lemma~\ref{Lem:3.3} that
$l_{\h}(\g/\h)\leqslant 1$ implies $i(\h,\g)=1$.

In Table~\ref{Tbl:3.4} we list all (up to the conjugacy by an
automorphism)
 simple subalgebras $\h$ in classical  simple algebras $\g$
with $i(\h,\g)=1$.

\begin{longtable}{|c|c|}
\caption{Simple subalgebras $\h$ in classical algebras $\g$ with
$\iota(\h,\g)=1$}\label{Tbl:3.4}\\\hline
$\g$&$\h$\\\endfirsthead\hline $\g$&$\h$\\\endhead\hline $\sl_n,
n\geqslant 2$&$\sl_k, k\leqslant n$\\\hline $\sl_n, n\geqslant
4$&$\sp_{2k}, 2\leqslant k\leqslant n/2$\\\hline $\so_n,n\geqslant
7$&$\so_k, k\leqslant n, k\neq 4$\\\hline $\so_n,n\geqslant
7$&$\sl_k, k\leqslant n/2$\\\hline $\so_n,n\geqslant 8$&$\sp_{2k},
2\leqslant k\leqslant n/4$\\\hline $\so_n,n\geqslant
7$&$G_2$\\\hline $\so_n, n\geqslant 9$&$\spin_7$\\\hline
$\sp_{2n},n\geqslant 2$&$\sp_{2k},k\leqslant n$\\\hline
\end{longtable}

\begin{Prop}\label{Prop:3.5}
Simple proper subalgebras $\h\subset\g$ (up to conjugacy in
$\Aut(\g)$) such that $l_\h(\g/\h)< 1$ (respectively,
$l_{\h}(\g/\h)=1$), are exactly those given in Table~\ref{Tbl:3.6}
(respectively, in Table~\ref{Tbl:3.7}).
\end{Prop}
\begin{proof}
This is checked directly using Table~\ref{Tbl:3.2}, equality
(\ref{eq:3.1}), Table~\ref{Tbl:3.4} (for classical $\g$) or results
of the classification in~\cite{Dynkin} (for  exceptional $\g$).
\end{proof}

 In the last column of Table~\ref{Tbl:3.7} the nontrivial part
of the $\h$-module $\g/\h$ is given. Here $\tau$  denotes the
tautological representation of a classical Lie algebra.

\begin{longtable}{|c|c|c|}
\caption{Simple subalgebras $\h\subsetneqq\g$ with $l_\h(\g/\h)<
1$}\label{Tbl:3.6}\\\hline N&$\g$&$\h$\\\endfirsthead\hline
N&$\g$&$\h$\\\endhead\hline 1&$\sl_n, n>1$&$\sl_k, n/2< k<
n$\\\hline 2&$\sl_{2n},n\geqslant 2$&$\sp_{2n}$\\\hline
3&$\sl_{2n+1},n\geqslant 2$&$\sp_{2n}$\\\hline
4&$\sp_{2n},n\geqslant 2$&$\sp_{2k}, n/2\leqslant k< n$\\\hline
5&$\so_n, n\geqslant 7$&$\so_k, \frac{n+2}{2}< k< n, k\neq
4$\\\hline 6&$\so_{2n},n\geqslant 5$&$\sl_n$\\\hline
7&$\so_{2n+1},n\geqslant 3 $&$\sl_n$\\\hline 8&$\so_n, 9\leqslant
n\leqslant 11$&$\so_7$\\\hline 9&$\so_n, 7\leqslant n\leqslant
9$&$G_2$\\\hline 10&$G_2$&$A_2$\\\hline 11&$F_4$&$B_4$\\\hline
12&$F_4$&$D_4$\\\hline 13&$F_4$&$B_3$\\\hline 14&$E_6$&$F_4$\\\hline
15&$E_6$&$D_5$\\\hline 16&$E_6$&$B_4$\\\hline 17&$E_7$&$E_6$\\\hline
18&$E_7$&$D_6$\\\hline 19&$E_8$&$E_7$\\\hline
\end{longtable}

\begin{longtable}{|c|c|c|c|}
\caption{Simple subalgebras $\h\subsetneqq\g$ with
$l_\h(\g/\h)=1$}\label{Tbl:3.7}\\\hline
N&$\g$&$\h$&$\g/\h_+$\endfirsthead\hline
N&$\g$&$\h$&$\g/\h_+$\\\endhead\hline 1&$\sl_{2n},
n>1$&$\sl_n$&$n(\tau+\tau^*)$\\\hline 2&$\sl_{2n},n\geqslant
3$&$\sp_{2n-2}$&$R(\pi_2)+4\tau$\\\hline 3&$\sp_{4n+2},n\geqslant
1$&$\sp_{2n}$&$(n+2)\tau$\\\hline 4&$\so_{2n}, n\geqslant
4$&$\so_{n+2}$&$(n-2)\tau$\\\hline
 5&$\so_{2n},n\geqslant
4$&$\sl_{n-1}$&$\bigwedge^2\tau+\bigwedge^2\tau^*+2(\tau+\tau^*)$\\\hline
6&$\so_{12}$&$\so_7$&$\tau+4 R(\pi_3)$\\\hline
7&$\so_{10}$&$G_2$&$4R(\pi_1)$\\\hline 8&$G_2$&$A_1$&$4\tau$\\\hline
9&$E_6$&$D_4$&$2(\tau+R(\pi_3)+R(\pi_4))$\\\hline
10&$E_6$&$A_5$&$2\bigwedge^3\tau$\\\hline
11&$E_7$&$F_4$&$3\tau$\\\hline
12&$E_7$&$B_5$&$\tau+2R(\pi_5)$\\\hline
\end{longtable}

Now we deduce two corollaries from Proposition~\ref{Prop:3.5}
\begin{Cor}\label{Cor:3.8}
Let $\g$ be a simple Lie algebra and $\h$ a semisimple subalgebra in
$\g$ such that $l_{\h_1}(\g/\h_1)<1$ for any simple ideal
$\h_1\subset \h$. If $\h$ is not simple, then $\g\cong\sp_{4n},
\h\cong\sp_{2n}\times\sp_{2n}$.
\end{Cor}
\begin{proof}
Checked immediately using Table~\ref{Tbl:3.6}.
\end{proof}

\begin{Cor}\label{Cor:3.9}
Let $\g$ be a simple Lie algebra and  $\h$ a semisimple subalgebra
of $\g$ such that $l_{\h_1}(\g/\h_1)=1$ for any simple ideal
$\h_1\subset \h$. The algebra $\h$ is not simple precisely in the
following three cases:
\begin{enumerate}
\item $\g\cong\sl_{2n}, n>1$, $\h\cong\sl_n\times\sl_n$.
\item $\g\cong\sp_{4n+2},$, $\h\cong \sp_{2n}\times\sp_{2n}$.
\item $\g\cong\sp_6$, $\h\cong\sl_2\times\sl_2\times\sl_2$.
\end{enumerate}
\end{Cor}
\begin{proof}
Checked immediately using Table~\ref{Tbl:3.7}.
\end{proof}

\section{Classification of initial subalgebras}
This section is central in the article. Here $\g$ is a semisimple
Lie algebra, $\h$ is an initial  (that is, saturated and
indecompsable) subalgebra of $\g$ and $\l_0$ is the s.s.g.p for the
$\h$-module $\g/\h$.

\begin{Prop}\label{Prop:4.1}
Precisely one of the following  possibilities {\rm (a),(b)} takes
place:
\begin{itemize}
\item[(a)]  $\g$ is simple, $\z(\h)=0$ and $l_{\h_i}(\g/\h_i)=1$
for any simple ideal $\h_i\subset\h$.
\item[(b)] Any simple ideal $\g_1\subset\g$ satisfies  precisely one of the following conditions:
\begin{itemize}
\item[(b1)] There exists a simple non-commutative ideal $\h_1\subset\h$
contained in $\g_1$ such that $l_{\h_1}(\g_1/\h_1)<1$.
\item[(b2)] There exists a simple non-commutative ideal  $\h_1\subset \h$ contained in the sum of $\g_1$ and another
simple ideal of $\g$  and projecting to $\g_1$ isomorphically.
\end{itemize}
\end{itemize}
\end{Prop}
\begin{proof}
Clearly, possibilities (a),(b) can not  take place simultaneously.
The same is true  for conditions (b1),(b2) (for a given ideal
$\g_1\subset \g$).

Suppose now that (b) does not take place. In this case
Corollary~\ref{Cor:2.13} implies there exists a simple ideal
$\g_1\subset\g$ such that

1) $\g_1$ does not contain a simple ideal $\h_1\subset\h$ with
$l_{\h_1}(\g_1/\h_1)<1$.

2) There exist no simple ideals  $\h_1\subset\g$ projecting
surjectively onto $\g_1$.

Let $\pi_2$ denote the orthogonal projection $\g\rightarrow
\g_1^{\perp}$. Put $\widetilde{\g}=\pi_2^{-1}(\pi_2(\h))$. Note that
$\h\subset\widetilde{\g}$. Denote by $\h_0$ the kernel of the
representation of $\h$ in $\widetilde{\g}/\h$. Clearly, $\h_0$ is an
ideal of $\widetilde{\g}$, $\h_0=\h\cap\g_1^\perp$. Put
$\widehat{\g}=\widetilde{\g}/\h_0, \widehat{\h}=\h/\h_0$. Since
$\widetilde{\g}/\h$ is an $\h$-submodule in $\g/\h$, it follows that
$\widehat{\h}$ is an essential subalgebra in $\widehat{\g}$ (by
Lemma~\ref{Lem:2.9}). Note that, by the construction, the projection
of any ideal of $\widehat{\h}$ to $\g_1$ does not coincide with
$\{0\}$, $\g_1$.

Suppose that the projection of a simple ideal
$\h_1\subset\widehat{\h}$ to some simple ideal $\g_2\subset
\widehat{\g}, \g_2\neq \g_1$ is nonzero. Then $\h_1$ is an
$\h_1$-submodule in $\widehat{\g}/\h_1$. If, in addition,
$\l_{\h_1}(\widehat{\g}/\h_1)\leqslant 1$, then the $\h_1$-module
$\widehat{\g}/\h_1$ is the direct sum of $\h_1$ and the trivial
module because $l_{\h_1}(\h_1)=1$. Thus
$\g_1=\pi_1(\h_1)+\z_{\g_1}(\pi_1(\h_1))$. By the last equality,
$\pi_1(\h_1)$ is an ideal of $\g_1$.  Contradiction with condition
2) above. Hence $l_{\h_1}(\widehat{\g}/\h_1)>1$.

Furthermore, condition  1) implies that for any simple
non-commutative ideal $\h_1\subset\h$ contained in $\g_1$ the
inequality $l_{\h_1}(\widehat{\g}/\h_1)\geqslant 1$ holds. It
follows from Proposition~\ref{Prop:3.1} that $\z(\widehat{\h})=0$,
and  $l_{\h_1}(\widehat{\g}/\h_1)=1$ for all simple ideals
$\h_1\subset\widehat{\h}$. We deduce that all simple ideals of
$\widehat{\h}$ are contained in $\g_1$. Thus $\widehat{\g}=\g_1$.
Since $\h$ is indecomposable,  possibility (a) takes place.
\end{proof}

All initial subalgebras $\h\subset \g$ satisfying condition  (a) of
Proposition~\ref{Prop:4.1} can be easily described: we classify them
in Proposition~\ref{Prop:4.5}. The classification of initial
subalgebras $\h$ satisfying condition (b) is more difficult. The
following proposition describes initial subalgebras possessing
ideals of some special form. In Proposition~\ref{Prop:4.6} we will
see that almost all simple ideals $\h_1\subset\h$ contained in some
simple ideal $\g_1\subset\g$ and satisfying the inequality
$l_{\h_1}(\g_1/\h_1)<1$ have this form.

\begin{Prop}\label{Prop:4.2}
Let $\h_1$ be a semisimple  ideal of $\h$ contained in an ideal
$\g_1\subset\g$. Put $\widehat{\h}_1=\n_{\g_1}(\h_1)$ and denote by
$\f$ the s.s.g.p. for the $\widehat{\h}_1$-module
$\g_1/\widehat{\h}_1$. Suppose $[\f,\f]\subset\h_1$. Then $\g=\g_1$
and $\h$ is an ideal in $\widehat{\h}_1$ generated by $\l_0$.
\end{Prop}
To prove this proposition we need the following lemmas.
\begin{Lem}\label{Lem:4.3}
Let $(\g,\h,\f)$ be a  triple of reductive Lie algebras such that
$\f\subset \h\subset \g$ and $\h$ contains the ideal of $\g$
generated by $[\f,\f]$. Denote by $\f_0$ the s.s.g.p. for the
$\f$-module $\g/\h$.  Then $\f_0$ is the intersection of $\f$ with
the maximal ideal of $\g$ contained in $\h$. In particular, $\h$
contains an ideal of $\g$ generated by  $\f_0$.
\end{Lem}
\begin{proof}
It is enough to prove the lemma in the case when $\h$ does not
contain an ideal of $\g$. In this case one should prove
$\f_0=\{0\}$. The algebra $\f$ is commutative. Thus $\f_0$ coincides
with the inefficiency kernel of the $\f$-module $\g/\h$.
 We see that $\g=\h+\z_\g(\f_0)$. Taking into account that $\h$ is
$\ad(\f_0)$-stable, we get $[\g,\f_0]\subset \h$. Since
$\g=\z_\g(\f_0)\oplus [\g,\f_0]$ and $[\z_\g(\f_0),[\g,\f_0]]\subset
[\g,\f_0]$, the subalgebra of $\g$ generated by $[\g,\f_0]$ is an
ideal. Hence $\h$ contains the ideal of $\g$ generated by
$[\g,\f_0]$. By our assumption, $\h$ does not contain an ideal of
$\g$. This implies $[\g,\f_0]=\{0\}$. Analogously,
$\h\cap\z(\g)=\{0\}$. So $\f_0\subset \z(\g)\cap \h=\{0\}$.
\end{proof}

\begin{Lem}\label{Lem:4.4}
For a saturated subalgebra $\h\subset \g$ the equality
$\h=[\h,\h]^{sat}$ holds.
\end{Lem}
\begin{proof}
By the second assertion  of Corollary~\ref{Cor:2.12},
$\h_0:=\h/[\h,\h]$ is an essential subalgebra of
$\g_0:=\n_\g([\h,\h])/[\h,\h]$. Since $\h_0$ is commutative, the
s.s.g.p. for  the $\h_0$-module $\g_0/\h_0$ is the inneficiency
kernel of this module. Thus $\h_0\subset \z(\g_0)$. Denote by
$\widetilde{\h}$ (resp., by $\widehat{\h}$) the inverse image of
$\z(\g_0)$ in $\n_\g([\h,\h])$ (resp., of $\z(\n_\g(\h)/\h)$ in
$\n_\g(\h)$) under the natural projection. By definition,
$[\h,\h]^{sat}=\widetilde{\h}^{ess}$. Note that $\z(\g_0)\subset
\z(\n_\g(\h)/[\h,\h])$. Thus $\widetilde{\h}\subset \widehat{\h}$.
The inclusions $\h\subset \widetilde{\h}\subset \widehat{\h}$ imply
$\h^{ess}\subset \widetilde{\h}^{ess}=[\h,\h]^{sat}\subset
\widehat{\h}^{ess}$ (see Corollary~\ref{Cor:2.7}). At last, we note
that $\h^{ess}=\h, \widehat{\h}^{ess}=\h$ because $\h$ is saturated.
\end{proof}

\begin{proof}[Proof of Proposition~\ref{Prop:4.2}]
Denote by $\pi_1,\pi_2$ the orthogonal projections $\g\rightarrow
\g_1,\g\rightarrow \g_1^\perp$.

{\it Step 1.} Let us prove that any simple non-commutative ideal
$\a\subset \h$ is contained either in $\g_1$ or in $\g_1^\perp$.
Assume the converse. Put $\g_0=\pi_2^{-1}(\pi_2(\a))$. This is a
semisimple Lie algebra isomorphic to $\g_1\oplus\a$. The subspace
$\h_0=\a\oplus \pi_1(\a^\perp\cap \h)\subset\g_0$ is a subalgebra
because $[\a,\pi_1(\a^\perp\cap \h)]=0$. There is the epimorphism
$\varphi:\h\rightarrow \h_0, \varphi|_{\a}=id,
\varphi|_{\a^\perp\cap\h}=\pi_1$. The $\h$-module $\g_0/\h_0$
 is isomorphic to
$\a\oplus \g_1/\pi_1(\h_1)$. Thus $\g_0/\h_0$ is a submodule in
$\g/\h$. It follows from assertion 1 of Corollary~\ref{Cor:2.12}
that $\h_0$ is an essential subalgebra of $\g_0$.

 Let
$\f_1$ denote the s.s.g.p. for the $\h_0$-module
$\g_1/\widehat{\h}_1\cong \g_0/\n_{\g_0}(\h_1)\hookrightarrow
\g_0/\h_0$. Now we show that the triple
$(\n_{\g_0}(\h_1),\h_0,\f_1)$ satisfies the assumptions of
Lemma~\ref{Lem:4.3}. By the construction, $\h_1$ is an ideal in
$\h_0$. Therefore $\h_0\subset \n_{\g_0}(\h_1)$ and $\h_1$ is an
ideal of $\n_{\g_0}(\h_1)$ contained in $\h_0$. The restriction of
the homomorphism $\pi_1:\g\rightarrow\g_1$ to $\h_0$ is an embedding
and $\pi_1(\h_0)=\pi_1(\h)$. The subalgebra $\pi_1(\f_1)\subset
\pi_1(\h)$ coincides with the s.s.g.p for the $\pi_1(\h)$-module
$\g_1/\widehat{\h}_1$.  It follows from Lemma~\ref{Lem:2.4} that
$\pi_1(\f_1)\subset \Ad(g)\f$ for some $g\in \Int(\widehat{\h}_1)$.
Thus $[\f_1,\f_1]\subset \h_1$. So we have checked that the triple
$(\n_{\g_0}(\h_1),\h_0,\f_1)$ satisfies the assumptions of
Lemma~\ref{Lem:4.3}.

Let $\f_0$ be the s.s.g.p. for the $\f_1$-module
$\n_{\g_0}(\h_1)/\h_0$.  By Lemma~\ref{Lem:2.9}, $\f_0$ is also the
s.s.g.p for the $\h_0$-module $\g_0/\h_0$. Let us show that the
projection of $\f_0$ to $\g_0\cap \g_1^\perp$ is nonzero. Assume the
converse. Then $\f_0\subset \g_1\cap\h_0=\h_0\cap \a^{\perp}$. Hence
the projection of $\f_0$ to $\a$ is zero. Since $\h_0$ is an
essential subalgebra of $\g_0$, we get the contradiction. But
$\g_1^\perp\cap \g_0$  is an ideal of $\n_{\g_0}(\h_1)$. By
Lemma~\ref{Lem:4.3} applied to $(\n_{\g_0}(\h_1),\h_0,\f_1)$, we
have $\g_1^\perp\cap \g_0\subset \h_0$ which is absurd.

{\it Step 2.} By step 1,  $[\h,\h]=([\h,\h]\cap \g_1^\perp)\cap
([\h,\h]\cap\g_1)$. By Lemma~\ref{Lem:4.4},
$\h=(\h\cap\g_1^\perp)\oplus (\h\cap \g_1)$. Since $\h$ is
indecomposable, $\g=\g_1$. Denote by $\f_1$ the s.s.g.p. for the
$\h$-module $\g/\widehat{\h}_1$. Note that $\f_1$ is
$\Int(\widehat{\h}_1)$-conjugate to a subalgebra of $\f$. Thus
$[\f_1,\f_1]\subset \h_1$ and  the triple $(\widehat{\h}_1,\h,\f_1)$
satisfies the assumptions of Lemma~\ref{Lem:4.3}. The s.s.g.p. for
the $\f_1$-module $\widehat{\h}_1/\h$ is $\Int (\h)$-conjugate to
$\l_0$. Applying Lemma~\ref{Lem:4.3}, we see that $\h$ contains the
ideal of $\widehat{\h}_1$ generated by $\l_0$. Since $\h$ is
essential, it coincides with this ideal.
\end{proof}

\begin{Prop}\label{Prop:4.5}
For the pair $(\h,\g)$ the following conditions are equivalent:
\begin{enumerate}
\item $(\h,\g)$ satisfies condition (a) of
Proposition~\ref{Prop:4.1} and $\h$ is an initial subalgebra of
$\g$.
\item $(\h,\g)=(\sl_n\times\sl_n,\sl_{2n}),(\sl_2\times\sl_2\times\sl_2,\sp_6), (\so_{n+2},\so_{2n}),
(A_5,E_6)$.
\end{enumerate}
\end{Prop}
\begin{proof}
Let $\h$ be an essential subalgebra of $\g$ satisfying condition (a)
of Proposition~\ref{Prop:4.1}.
 Note that $\h^{sat}$ also satisfies this condition. Therefore $\h^{sat}$ is semisimple.
Thus $\h^{sat}=\h$ and $\h$ is initial.

First  suppose $\h$ is not simple. Then $(\h,\g)$ is one of the
pairs listed in Corollary~\ref{Cor:3.9}.

Let $\tau_1,\ldots,\tau_i$ denote the tautological representations
of simple ideals of $\h$ ($i=2$ in cases 1,2 of
Corollary~\ref{Cor:3.9}, and  $i=3$ in case 3).

When $\h=\sl_n\times\sl_n, \g=\sl_{2n}$, the non-trivial part of the
representation of $\h$ in $\g/\h$ is
$\tau_1\otimes\tau_2+\tau_1^*\otimes\tau_2^*$. The s.s.g.p. for
$\tau_1\otimes\tau_2$ is the subalgebra $\sl_n\subset \h$ embedded
diagonally. The restriction of $\tau_1^*\otimes\tau_2^*$ to this
subalgebra is isomorphic to the direct sum of the adjoint and the
trivial representations. Thus the s.s.g.p. for the $\h$-module
$\g/\h$ is non-trivial and its projections to the both simple ideals
of  $\h$ are nonzero. By Proposition~\ref{Prop:2.6}, $\h$ is
essential.

Now let $\g=\sp_6, \h=\sl_2\times\sl_2\times\sl_2$. The
representation of  $\h$ in $\g/\h$ is isomorphic to
$\tau_1\otimes\tau_2+\tau_2\otimes\tau_3+\tau_3\otimes\tau_1$. The
s.s.g.p. $\l_1$ of the representation $\tau_2\otimes\tau_3$ is
isomorphic to $\sl_2\times\sl_2$, where the second factor is
embedded diagonally into the sum of two ideals $\sl_2\subset \h$.
The restriction of the representation of $\h$ in $\g/\h$ to $\l_1$
is isomorphic to $2\tau_1\otimes\tau_2$. By the calculation in the
previous case, the s.s.g.p. for this representation is
one-dimensional and projects non-trivially to the both simple ideals
of $\l_1$. Thus $\h$ is essential.

In the remaining case we have to prove that the s.s.g.p for the
$\h$-module $\g/\h$ is trivial. Let $\g=\sp_{4n+2},
\h=\sp_{2n}\times \sp_{2n}$. The non-trivial part of the
representation of $\h$ in $\g/\h$ is
$\tau_1\otimes\tau_2+2(\tau_1+\tau_2)$. The s.s.g.p. $\l_1$ for
$\tau_1\otimes\tau_2$ is the direct sum of  $n$ copies of $\sl_2$,
embedded diagonally to $\sp_{2n}\times\sp_{2n}$ (see~\cite{Elash2}).
It can be easily seen that the restriction of
$2(\tau_1\otimes\tau_2)$ to $\l_1$ has the trivial s.s.g.p.

If $\h$ is simple, then the proof follows from Table~\ref{Tbl:3.7}
and the classification in~\cite{Elash1}.
\end{proof}

Until the end of the section  we suppose that  $\g$ is a semisimple
Lie algebra and $\h$ its initial subalgebra satisfying condition (b)
of Proposition~\ref{Prop:4.1}. The latter means that for any simple
ideal  $\g_1\subset \g$ there exists a simple non-commutative ideal
$\h_1\subset \h$ either contained in $\g_1$ and satisfying the
inequality $l_{\h_1}(\g_1/\h_1)<1$ or projecting to  $\g_1$
isomorphically  and contained in the sum of $\g_1$ and another
simple ideal of $\g$.

\begin{Prop}\label{Prop:4.6}
Let  $\g_1$ be a simple ideal of $\g$ and $\h_1$ a simple
non-commutative ideal of $\h$ such that $\h_1\subset \g_1$ and
$l_{\h_1}(\g_1/\h_1)<1$. If $\g_1,\h_1$ do not satisfy the
assumption of Proposition~\ref{Prop:4.2}, then $\g_1\cong \sp_{2m},
m\geqslant 2$.
\end{Prop}
\begin{proof}
Put $\widehat{\h}_1=\n_{\g_1}(\h_1)$ and denote by $\l_1$ the
s.s.g.p. for the  $\widehat{\h}_1$-module $\g_1/\widehat{\h}_1$.

Assume the converse, let $\g_1\not\cong\sp_{2m}$ and $[\l_1,\l_1]$
is not contained in $\h_1$. By Corollary~\ref{Cor:3.8},  if $\h_2$
is an ideal in $\widehat{\h}_1$ such that $l_{\g_1}(\g_1/\h_2)<1$,
then $\h_2=\h_1$. Thus Lemma~\ref{Lem:2.4} and assertion 1 of
Proposition ~\ref{Prop:3.1} imply that $[\l_1,\l_1]$ projects
injectively into $\h_1$.

By Corollary~\ref{Cor:2.7}, $[\l_1,\l_1]$ coincides with the
commutant of some Levi subalgebra of $\g_1$. Fix a simple ideal
$\f\subset [\l_1,\l_1]$ and let $\h_1,\ldots,\h_k$ be all simple
ideals of $\widehat{\h}_1$ such that projection of $\f$ to $\h_i$ is
nonzero, or, equivalently, injective. Then $i(\f,\g)=\sum_j
i(\f,\h_j)i(\h_j,\g)$, see \cite{Dynkin}. Hence for any simple ideal
$\f\subset [\l_1,\l_1]$ not contained in $\h_1$ the inequality
$i(\f,\g)\geqslant 2$ holds. This implies $\g_1\not\cong
\sl_n,\so_{2n},E_l$ because the Dynkin index of any simple ideal of
a Levi subalgebra in $\sl_n,\so_{2n},E_l$ equals 1. If  $\g_1\cong
F_4,G_2$, then Table~\ref{Tbl:4.8} implies
$[\widehat{\h}_1,\widehat{\h}_1]=\h_1$. It remains to consider the
case $\g\cong\so_{2n+1},n\geqslant 3$. We may assume that
$[\widehat{\h}_1,\widehat{\h}_1]$ is not simple. So we have to
consider only cases 5, 8 (n=11).

Consider case  5 from Table~\ref{Tbl:4.8}. In this case the
$\widehat{\h}_1$-module $\g/\widehat{\h}_1$ is isomorphic to the
tensor product of the tautological  $\so_k$ and
$\so_{2n+1-k}$-modules. By~\cite{Elash2}, $\l_1\subset \so_k$.

In case 8 the algebra $\widehat{\h}_1$ has two simple ideals and is
contained (up to conjugacy) in $\so_8\times \so_3$. By the previous
paragraph, the projection of $\l_1$ to the ideal $\so_3\subset
\widehat{\h}_1$ is zero.
\end{proof}

Recall that the {\it complexity} of an irreducible $G$-variety $X$
is the codimension of a general $B$-orbit, or, equivalently, $\td
K(X)^B$. We denote the complexity of $X$ by  $c_G(X)$. A normal
$G$-variety $X$ of complexity 0 is called {\it spherical}. An
algebraic  subalgebra $\h\subset \g$ is said to be spherical if
$G/H$ is spherical.

In the sequel we need the following standard
\begin{Prop}\label{Prop:4.7}
Let $H\subset \widetilde{H}\subset G$ be connected reductive
algebraic groups such that $\widetilde{H}\subset N_G(H)$ and
$\widetilde{H}/H$ is a torus. Note that $\widetilde{H}/H$ acts on
$G/H$ by $G$-automorphisms. Then the representation of
$\widetilde{H}/H$ in $K[G/H]^U$ is effective. Further, the following
conditions are equivalent
\begin{enumerate}
\item $\rank_G(G/H)=\rank_G(G/\widetilde{H})+\dim \widetilde{H}/H$.
\item $c_G(G/H)=c_{G}(G/\widetilde{H})$
\item The action of $\widetilde{H}/H$ on $K(G/H)^B$ is trivial.
\item For a highest weight $\lambda$ the group $\widetilde{H}/H$
acts on $K[G/H]^U_{(\lambda)}$ by the multiplication by a character.
\item For a highest weight $\lambda$ the group $\widetilde{H}/H$
acts on $V(\lambda)^{H}$ by the multiplication  by a character.
\end{enumerate}
\end{Prop}
\begin{proof}
If $h\in \widetilde{H}$ acts trivially on $K[G/H]^U$, then it acts
trivially on $K[G/H]$ because $h$ acts as a $G$-automorphism and the
$G$-module $K[G/H]$ is generated by $K[G/H]^U$. This implies $h\in
H$.

The equivalence $(2)\Leftrightarrow(3)$ is straightforward.
 The equivalence of (4)
and (5) follows from the Frobenius reciprocity. Implication
$(4)\Rightarrow (3)$ is tautological. By \cite{VP}, Theorem 3.14,
$K[G/H]^U_{\lambda}$ is a finite-dimensional vector space. So it is
spanned by $\widetilde{H}/H$-eigenvectors. To prove $(3)\Rightarrow
(4)$ it remains to note that $f_1/f_2\in K(G/H)^B$ for any
$f_1,f_2\in K[G/H]^{U}_{(\lambda)}$.

It remains to prove the equivalence of (1) and (2).

By the results of~\cite{Knop1}, Abschnitt 7,
\begin{equation}\label{eq:compl_formula}
\begin{split}
&c_G(G/H)=\frac{1}{2}(2(\dim G-\dim H)-(\dim G-\dim L_{0\,G,G/H})-\rank_G(G/H))=\\
&=\frac{1}{2}(\dim G+\dim L_{G,G/H})-\dim H-\rank_G(G/H).
\end{split}
\end{equation}
The commutants of $\l_{G,G/H},$ $\l_{G,G/\widetilde{H}}$ are
conjugate to the commutants of the s.s.g.p's for the $\h$-module
$\g/\h$ and the $\widetilde{\h}$-module $\g/\widetilde{\h}$,
respectively. Since the last two commutants coincide, we see that
$L_{G,G/H}=L_{G,G/\widetilde{H}}$ and
$c_G(G/H)-c_G(G/\widetilde{H})=\dim\widetilde{H}/H+\rank_G(G/\widetilde{H})-\rank_G(G/H)$.
Thus $(1)\Leftrightarrow (2)$ is proved.
\end{proof}

In the following table we give the normalizers  $\n_\g(\h_1)$ of
subalgebras  $\h_1\subset\g$ listed in Table~\ref{Tbl:3.6} (in the
column $\h$), the  $\n_\g(\h_1)$-modules $\g/\n_\g(\h_1)$ and the
ideals of $\n_\g(\h_1)$ generated by their intersections with
$\l_1$, where
 $\l_1$ denotes the s.s.g.p for the $\n_\g(\h_1)$-module $\g/\n_\g(\h_1)$.

In the first and the second columns the number of the corresponding
subalgebra $\h_1$ in Table~\ref{Tbl:3.6} and $\n_\g(\h_1)$ are
given. The first ideal of $\n_\g(\h_1)$ in the second  column is
$\h_1$. In all cases  $\dim \z(\n_\g(\h_1))\leqslant 1$. If the
equality holds, then, conjugating $\h_1$ in $\g$ if necessary, we
may assume that $\z(\n_\g(\h_1))$ lies in $\t$ and is spanned by the
dual fundamental weight.

The representation of $\n_\g(\h_1)$ is written as the sum of
irreducible subrepresentations. The order of irreducible
representations over simple algebras in tensor products coincides
with the order of simple ideals of $\n_\g(\h_1)$ in column 2. The
lower index indicates the constant, by that the element of
$\z(\n_\g(\h_1))$ given in column  2 in the triangular brackets acts
on the corresponding irreducible module. By $\tau_1,\ldots,\tau_k$
we denote the tautological representations of classical simple
ideals of $[\n_\g(\h_1),\n_\g(\h_1)]$, and by ${\bf 1}$ the
one-dimensional trivial representation.

Ideals of $\n_\g(\h_1)$ generated by their intersections with $\l_1$
are listed in column 4. In all  cases, except N4, any initial
subalgebra $\h\subset\g$ such that  $\h_1$ is an ideal of $\h$ is
the one of the listed ideals
(Propositions~\ref{Prop:4.2},\ref{Prop:4.6}). Such an ideal is given
by a sequence of numbers. The sequences are separated by "$;$". The
sequence $i_1,\ldots,i_j$ corresponds to the sum of ideals with
numbers $i_1,\ldots,i_j$. Here the number of an ideal is its
position in column 2. For example, for
$\n_\g(\h_1)=\sl_k\times\sl_{n-k}\times\langle\pi_k^\vee\rangle$ the
ideal corresponding to the sequence $"1,3"$ is
$\sl_k\times\langle\pi_k^\vee\rangle$.

\begin{longtable}{|c|c|c|c|}
\caption{The representations of $\n_\g(\h_1)$ in
$\g/\n_\g(\h_1)$}\label{Tbl:4.8}\\\hline
N&$\n_\g(\h_1)$&$\g/\n_\g(\h_1)$&$\h$\\\endfirsthead\hline
N&$\n_\g(\h_1)$&$\g/\n_\g(\h_1)$&$\h$\\\endhead\hline
1&$\sl_k\times\sl_{n-k}\times\langle\pi_k^\vee\rangle$&$_1\tau_1\otimes\tau_2^*+_{-1}\tau_1^*\otimes\tau_2$&$1(k\neq
\frac{n+1}{2}) ;1,2;$\\&&&$1,3(k\neq n/2);1,2,3(k\neq n/2)$\\\hline
2&$\sp_{2n}$&$R(\pi_2)$&1\\\hline
3&$\sp_{2n}\times\langle\pi_{2n}^\vee\rangle$&$R(\pi_2)+_1\tau+_{-1}\tau$&1,2\\\hline
4&$\sp_{2k}\times\sp_{2n-2k}$&$\tau_1\otimes\tau_2$&1 ($n\neq 2k
$);1,2\\\hline
5&$\so_k\times\so_{n-k}$&$\tau_1\otimes\tau_2$&1\\\hline
6&$\sl_n\times
\langle\pi_{n}^{\vee}\rangle$&$_1\bigwedge^2\tau+_{-1}\bigwedge^2\tau^*$&1;1,2
$(n=2m+1)$
\\\hline
7&$\sl_n\times\langle\pi_n^\vee\rangle$&$_2\bigwedge^2\tau+_{-2}\bigwedge^2\tau+_1\tau+_{-1}\tau^*$&\\\hline
8&$\so_7\times\so_{n-8}$&$\tau_1+R(\pi_3)\otimes\tau_2$&1
($n\leqslant 10$)\\\hline
9&$G_2\times\so_{n-7}$&$R(\pi_1)\otimes(\tau_2+{\bf 1})$&1
($n\leqslant 8$)
\\\hline 10&$A_2$&$\tau+\tau^*$&1\\\hline
11&$B_4$&$R(\pi_4)$&1\\\hline
12&$D_4$&$\tau+R(\pi_3)+R(\pi_4)$&1\\\hline
13&$B_3\times\langle\pi_1^\vee\rangle$&$_1\tau+_{-1}\tau+_1R(\pi_3)+_{-1}R(\pi_3)$&\\\hline
14&$F_4$&$R(\pi_1)$&1\\\hline
15&$D_5\times\langle\pi_1^{\vee}\rangle$&$_1R(\pi_4)+_{-1}R(\pi_5)$&1;1,2\\\hline
16&$B_4\times \langle
\pi_1^\vee\rangle$&$_1R(\pi_4)+_{-1}R(\pi_4)+\tau$&1\\\hline
17&$E_6\times\langle\pi_1^{\vee}\rangle$&$_1R(\pi_1)+_{-1}R(\pi_5)$&1\\\hline
18&$D_6\times A_1$&$R(\pi_5)\otimes\tau$&1\\\hline 19&$E_7\times A_1
$&$R(\pi_1)\otimes\tau$&1\\\hline
\end{longtable}
\begin{Prop}\label{Prop:4.9}
Ideals  $\h\subset\n_\g(\h_1)$ containing $\h_1$ and generated by
$\l_1\cap\h$ are precisely those listed in Table~\ref{Tbl:4.8}.
\end{Prop}
\begin{proof}
If $\n_\g(\h_1)$ is simple or the $\n_\g(\h_1)$-module
$\g/\n_\g(\h_1)$ is irreducible we use the results
of~\cite{Elash1},~\cite{Elash2}, respectively. Note that in case 17
there is a mistake in~\cite{Elash2}:  the s.s.g.p. is contained in
$D_6$. Let us prove this.

Consider a linear action of a reductive group $H_1\times H_2$, where
$H_1$ is semisimple,  on a vector space $V$. The projection of the
s.g.p. for this action to $H_2$ is the s.g.p. for the action
$H_2:V\quo H_1$. This is a consequence of the following fact: a
general fiber of the quotient morphism of the action  $H_1:V$
contains a dense orbit (see, for example,~\cite{Kraft}, ch.2, $\S$4,
D).

Now we note that $Y=(\g/\n_\g(\h_1))\quo\Spin(12)$ is an affine
space, the action of $\SL(2)$ on $Y$ is linear and the
$\SL(2)$-module $Y$ is the direct sum of the trivial two-dimensional
and the irreducible 5-dimensional $\SL(2)$-modules.  This follows
from the results of~\cite{Schwarz}. It remains to note that s.s.g.p.
for the $\SL(2)$-modul $Y$ is trivial.

It remains to consider  cases 1,3,6,7,8 for $n\geqslant 10$, 9 for
$n\geqslant 9$,13,15,16.

In cases 3,6,7,15,16 the required statement follows from the
classification in~\cite{Kramer}. Indeed, in all these cases the
subalgebra $\n_\g(\h_1)\subset\g$ is spherical, its center is
one-dimensional and the commutant is simple. By
Proposition~\ref{Prop:4.7}, the subalgebra $\h_1\subset \g$ is
spherical iff $\dim \a(\g,\h_1)>\dim \a(\g,\n_\g(\h_1))$ iff (see
Corollary~\ref{Cor:2.7}) the s.s.g.p's  for the $\h_1$-module
$\g/\h_1$ and the  $\n_\g(\h_1)$-module $\g/\n_\g(\h_1)$ are
distinct. We obtain the required list of ideals using Tabelle 1
from~\cite{Kramer}.

In the sequel we denote by $H_1$ the connected subgroup of $G$
corresponding to $\h_1$.

Consider case 1. By the main theorem from~\cite{Kraft}, ch.2, $\S4$,
for $k>n/2$ the action $N_G(H_1)/H_1:(\g/\n_\g(\h_1))\quo H_1$ is
the adjoint representation of $\GL_{n-k}$. Since the projection of
$\l_1$ to $\n_\g(\h_1)$ is the s.s.g.p. for the action
$N_G(H_1)/H_1:(\g/\n_\g(\h_1))\quo H_1$, we are done. If $k=n/2$,
then the $N_G(H_1)/H_1$-variety $(\g/\n_\g(\h_1))\quo H_1$ is the
direct sum of the adjoint module of $\GL_{n-k}$ and of a non-trivial
two-dimensional $\GL_{n-k}/\SL_{n-k}$-module. Thus the s.s.g.p. for
the $\n_\g(\h_1)$-module $\g/\n_\g(\h_1)$ is contained in
$[\n_\g(\h_1),\n_\g(\h_1)]$ whence the desired result.

In case 13 the action of the torus $(N_G(H_1)/H_1)^{\circ}$ on
$\g/\n_\g(\h_1)\quo H_1$ is non-trivial, thus $\l_0\subset\h_1$. By
the results from~\cite{Elash1}, $\l_0=0$.

 Analogously
we see that $\l_0=\{0\}$ in case 9 for $n=9$.

It remains to consider case 8. For $n=11$ the s.s.g.p. of the
representation $R(\pi_3)\otimes \tau_2$ of $\n_\g(\h_1)$ is
contained in $\h_1$ (\cite{Elash2}). The equality $\l_0=\{0\}$
follows from \cite{Elash1}. For $n=10$ the projection of $\l_0$ to
$\n_\g(\h_1)/\h_1$ is zero because the one-dimensional torus
$(N_G(H_1)/H_1)^\circ$ acts non-trivially on $(\g/\n_\g(\h_1))\quo
H_1$.
\end{proof}

Now it remains to find all initial subalgebras $\h\subset\g$, that
have no simple ideal $\h_1$ satisfying the conditions of
Proposition~\ref{Prop:4.2} with simple  $\g_1$. It follows from
Proposition~\ref{Prop:4.6} that any simple ideal of $\g$ satisfying
condition (b1) of Proposition~\ref{Prop:4.1} is isomorphic to
$\sp_{2n},n\geqslant 2$.

 First we consider the case $\g=\sp_{2n}$.
We classify initial subalgebras $\h\subset\sp_{2n}$, containing the
ideal $\h_1=\sp_{2k}, k\geqslant n/2$. We assume that $\h\neq
\h_1,\n_\g(\h_1)$ (these two subalgebras are initial, except the
case $\h=\h_1, k=n/2$; this follows from
Proposition~\ref{Prop:4.9}). Note that $\n_{\g}(\h_1)\cong
\sp_{2k}\times \sp_{2(n-k)}$ and that the $\n_{\g}(\h_1)$-module
$\g/\n_{\g}(\h_1)$ is the tensor product of the tautological
$\sp_{2k}$- and $\sp_{2(n-k)}$-modules.

The subalgebra $\h/\h_1\subset\n_\g(\h_1)/\h_1\cong\sp_{2(n-k)}$ is
essential (the third assertion of Corollary~\ref{Cor:2.12}). Thus it
contains an ideal $\h_2$ isomorphic to $\sp_{2l}, l\geqslant
(n-k+1)/2$ and so is contained in
$\sp_{2k}\times\sp_{2l}\times\sp_{2(n-k-l)}$.
\begin{Lem}\label{Lem:4.10}
Let $\h=\sp_{2k}\times\sp_{2l}\times\sp_{2(n-k-l)}\subset \g$ and
$\l_0$ be the s.s.g.p. for the $\h$-module $\g/\h$. Then the
following conditions are equivalent. \begin{enumerate}
\item $\l_0\not\subset\h_1$, \item $\h$ is essential, \item
$l=n-k-l=1$.\end{enumerate} Under these conditions  the projection
of the s.s.g.p. to $\h/\h_1$ is one-dimensional.
\end{Lem}
\begin{proof}
Let $\tau_1,\tau_2$ denote the tautological representations of
$\sp_{2l},\sp_{2(n-k-l)}$, respectively. By Theorem 7
from~\cite{Elash2}, the projection of  the s.s.g.p. for the
$\h$-module $\g/\n_\g(\h_1)$ to  $\z_\h(\h_1)\cong \sp_{2l}\times
\sp_{2(n-k-l)}$ is the s.s.g.p. for the representation
$\bigwedge^2(\tau_1+\tau_2)=\tau_1\otimes\tau_2+\bigwedge^2\tau_1+\bigwedge^2\tau_2$
of $\sp_{2l}\times\sp_{2(n-k-l)}$. Thus the projection of the
s.s.g.p. for the $\h$-module $\g/\h$ to $\z_\h(\h_1)$ coincides with
the s.s.g.p. for the representation $2\tau_1\otimes \tau_2+
\bigwedge^2\tau_1+\bigwedge^2\tau_2$ of $\z_\h(\h_1)$. The s.s.g.p.
for $\bigwedge^2\tau_1+\bigwedge^2\tau_2$ is $\l_1=(n-k)\sl_2\subset
\sp_{2l}\times\sp_{2(n-k-l)}$.  The index of the restriction of
$2\tau_1\otimes \tau_2$ to any simple ideal of $\l_1$ is not less
than 1. If $\s$ is a simple ideal of $\l_1$ contained in $\sp_{2l}$,
then this index is equal to $1$ iff $n-k-l=1$. To show that
$l,n-k-l\leqslant 1$ it is enough to use assertion 2 of
Proposition~\ref{Prop:3.1} and the fact that the s.s.g.p. for four
copies of the tautological $\sl_2$-module is trivial. Other
assertions are now clear.
\end{proof}
\begin{Cor}\label{Cor:4.11}
There exists a unique  subalgebra $\h\subset\n_\g(\h_1)$ essential
in $\g$ different from $\h_1,\n_\g(\h_1)$. This subalgebra is
$\sp_{2n-4}\times\sl_2\times\sl_2$.
\end{Cor}
\begin{proof}
By the previous lemma, we may assume $\h\subset
\widetilde{\h}:=\sp_{2n-4}\times\sl_2\times\sl_2$. The assertion
follows from Lemma~\ref{Lem:4.3} applied to the triple
$(\widetilde{\h},\h,\f)$, where $\f$ is the s.s.g.p. for the
$\widetilde{\h}$-module $\g/\widetilde{\h}$.
\end{proof}

Proceed to the general case. It remains to classify all initial
subalgebras $\h\subset \g$ such that $\g$ is not simple and for any
simple ideal $\g_1\subset\g$ exactly one of the following
possibilities takes place:

1) There exist  simple ideals $\h_1\subset\h,\g_2\subset\g,
\g_2\neq\g_1$ such that $\h_1\subset\g_1\oplus\g_2$ and $\h_1$
projects isomorphically onto $\g_1$.

2) $\g_1\cong\sp_{2n}$ and there exists an ideal
$\h_1\cong\sp_{2k}\subset\h, k\geqslant n/2$.

If there is no simple ideal in $\g$ satisfying condition 2), then
$\g\cong \g_1\times\g_1, \h\cong\g_1$ and $\h$ is embedded into $\g$
diagonally. This follows from the indecomposability of $\h$.

\begin{Lem}\label{Lem:4.12}
Suppose $\h$ is an initial subalgebra of $\g$ and   there exists an
ideal $\g_1^1\subset\g$ satisfying  condition 2). If $\g\not\cong
\sp_{2n}$, then $\g\cong \sp_{2m}\times \sp_{2n},
\h=\sp_{2m-2}\times\sl_2\times \sp_{2n-2}, n>1,m\geqslant 1$.
\end{Lem}
\begin{proof}
Denote by $\h_1^1$ a unique ideal of $\h$ embedded into $\g_1^1$ as
$\sp_{2n_1}, 2n_1\geqslant n$. Let $\a$ be a simple ideal of $\h$
not contained neither in $\g_1$ nor in $\g_1^\perp$. The existence
of such an ideal follows from Lemma~\ref{Lem:4.4} . By
Corollary~\ref{Cor:2.13}, there exists a simple ideal $\g^2_1\subset
\g$ such that $\a\subset \g_1:=\g^1_1\oplus \g_1^2$.

  If for $\g_1^2$ the first possibility listed above takes
place, then put $\h_1=\h_1^1\oplus \a$, $\h^2_1=\{0\}$. Otherwise,
$\g^2_1\cong \sp_{2m}$ and there is an ideal $\h_1^2\subset\h$
embedded into $\g_1^2$ as $\sp_{2m_1},2m_1\geqslant m$. We put
$\h_1=\h^1_1\oplus \h^2_1\oplus\a$.  Let us show that the pair
$(\g_1,\h_1)$ satisfies the conditions of
Proposition~\ref{Prop:4.2}. Let $\widehat{\h}_1=\n_{\g_1}(\h_1)$,
$\f$ be the s.s.g.p. for the $\widehat{\h}_1$-module
$\g_1/\widehat{\h}_1$ and   $\pi_i,i=1,2,$ denote the projection
$\g_1\rightarrow \g_1^i$. Let us check that
$\pi_1([\f,\f])\subset\pi_1(\a)\oplus \h_1^1$. Let
$\widehat{\h}^i_1=\n_{\g_1^i}(\pi_1(\a)\oplus \h^i_1), i=1,2$,
$\f^1$ be the s.s.g.p. for the $\widehat{\h}^1_1$-module
$\g_1^1/\widehat{\h}^1_1$. Recall that  the projection of $\h$ to
$\g_1^1\cong \sp_{2n}$ is an essential subalgebra in $\g_1$ (see
Corollary~\ref{Cor:2.12}). Inclusion $[\f^1,\f^1]\subset
\pi_1(\a)\oplus \h_1^1$ follows from the classification of essential
subalgebras in $\sp_{2n}$ obtained earlier. To prove the inclusion
$\pi_1([\f,\f])\subset \pi_1(\a)\oplus\h_1^1$ it remains to note
that $\pi_1(\f)$ is $\Int(\g)$-conjugate to a subalgebra $ \f^1$.
The latter follows from Lemma~\ref{Lem:2.9}. The inclusion
$\pi_2([\f,\f])\subset \pi_2(\a)\oplus\h^2_1$ is obvious (if
$\h^2_1=\{0\}$) or is proved analogously.   So we see that
$\pi_i([\f,\f])\subset \pi_i(\a)\oplus \h_1^i$. Thus $[\f,\f]\subset
(\h_1^1\oplus \pi_1(\a)\oplus\pi_2(\a)\oplus \h^1_2)\cap
\widehat{\h}_1=\h_1$.

Using Proposition~\ref{Prop:4.2},  we conclude that $\g=\g_1$ and
$\h$ is an ideal of $\widehat{\h}_1$ generated by $\f\cap \h$. In
particular, the projection of $\f$ to $\a$ is non-trivial.

Let us show that the projection $\f^1\rightarrow \pi_1(\a)$ is
surjective. Let $\a_0$ be the image of this projection. There is the
$\widehat{\h}_1$-submodule $V:=\a\oplus \g_1^1/\widehat{\h}^1_1$ in
$\g_1/\widehat{\h}_1$. For a Cartan subalgebra $\t\subset \a$ in
general position the intersection $\t\cap \a_0$ is the projection of
the s.s.g.p. for the $\widehat{\h}_1$-module $V$ to $\a$.  Let
$T,A,A_0$ denote the corresponding connected subgroups of $G$. Since
$\t$ is commutative, $\t\cap\a_0$ is the Lie algebra of the
inneficiency kernel $T_0$ of the action $T:A/A_0$. So $T_0\subset
gA_0g^{-1}$ for any $g\in A$. If $\t_0\neq 0$, then $A_0$ is a
normal subgroup in $A$. But the projection of $\f$ to $\a$ is
non-trivial, thus $\t\cap\a_0\neq 0$. We conclude that $\a_0=\a$.

Since the projection of $\f^1$ to $\pi_1(\a)$ is surjective,
$\pi_1(\a)\times \h^{1}_1\cong \sl_2\times \sp_{2n-2}$ (by our
classification of essential subalgebras in $\sp_{2n}$ and results
of~\cite{Elash2}).

Hence $\g_1^2\cong \sp_{2m}, m\geqslant 1$, $\pi_2(\a)\times
\h_1^2\cong \sl_2\times \sp_{2m-2}$.
\end{proof}

So a subalgebra $\h\subset \g$ is initial iff
$(\g,\h)=(\g_i,\h_i\oplus\z_i)$, where $(\g_i,\h_i)$ is listed in
Table~\ref{Tbl:1.6}, or  $(\g,\h)$ is a pair from
Table~\ref{Tbl:1.4} different from NN1,($k\neq n/2$),2,10,19.

\section{Computation of Cartan spaces and the classification of essential subalgebras}\label{SECTION_computation}

At first, we compute the spaces $\a(\g,\h)$ for semisimple essential
subalgebras.  It follows from Proposition~\ref{Prop:4.9} that all
semisimple essential indecomposable subalgebras are listed in
Table~\ref{Tbl:1.4}.

Firstly, we consider the case when $\g$ is  simple. When $\h$ is the
essential part of a spherical subalgebra of $\g$, the spaces
$\a_{\g,\h}$  can be computed using Tabelle 1 in~\cite{Kramer}.
There are only five pairs not satisfying this condition:
\begin{enumerate}
\item $\sl_2\times\sl_2\times\sl_2\subset\sp_6$.
\item $\sl_{k}\subset\sl_n, k\geqslant \frac{n+2}{2}$.
\item $\sp_{2n-4}\times\sl_2\times\sl_2\subset\sp_{2n}, n>3$.
\item $D_4\subset F_4$.
\item $B_4\subset E_6$
\end{enumerate}

In cases 1,3,4 the corresponding homogeneous space has complexity 1.
The spaces $\a(\g,\h)$ can be computed using results
of~\cite{Panyushev4}.

To compute $\a(\g,\h)$ in the remaining case we note that, by the
Frobenius reciprocity, $\a(\g,\h)$ is generated by highest weights
$\lambda$ such that $V(\lambda)^{*\h}\neq 0$. Since $\h$ is
reductive, the spaces $V(\lambda)^{\h}, V(\lambda)^{*\h}$ are dual
to each other. Therefore to compute $\a_{\g,\h}$ it is enough to
find $k=\dim\a_{\g,\h}$ independent highest weights
$\lambda_1,\ldots,\lambda_k$ such that $V(\lambda_i)^{\h}\neq 0,
i=\overline{1,k}$. Note that $\dim\a(\g,\h)=\rank\g-\rank\l_0$,
where $\l_0$ is the s.s.g.p for the $\h$-module $\g/\h$.

In case 2 $\dim\a(\g,\h)=2(n-k)$ because $\l_0\cong \sl_{2k-n}$.
Note that the highest vectors (for an appropriate choice of a Borel
subalgebra) in $V(\pi_i)$  for $i\leqslant n-k$ or $i\geqslant k$
are $\h$-invariant.

In case 5 $\dim V(\pi_1)^{\h}=2$. Since $V(\pi_2)=\bigwedge^2
V(\pi_1)$, it follows that $\pi_1,\pi_2\in \a(\g,\h)$. Further,
$\pi_5,\pi_4\in \a(\g,\h)$ because $V(\pi_5)=V(\pi_1)^*,
V(\pi_4)=\bigwedge^2 V(\pi_2)$. It is clear that $V(\pi_6)=\g$ has
an $\h$-invariant non-zero vector. Since $\l_0\cong\sl_2$, we see
that $\a(\g,\h)$ is generated by $\pi_1,\pi_2,\pi_4,\pi_5,\pi_6$.

It remains to consider the pairs NN25-27. In case 25 $\dim
V(\lambda)^\h\neq 0$ iff $V(\lambda)$ is the tensor product of two
$\h$-modules  dual to each other.

Now we consider  pairs 26,27. In this case $\l_0\cong
\sp_{2n-4}\times\sp_{2m-4}\times K$ and $\dim \a(\g,\h)=3$. It
remains to note that the $\g$-modules
$V(\pi_2),V(\pi_2'),V(\pi_1)\otimes V(\pi_1')$ contain non-zero
$\h$-invariant vectors.

This completes the proof of assertion (b) of Theorem~\ref{Thm:1.3}.
Assertion (d) is clear from Remark~\ref{Rem:1.7}.

Now we proceed to the proof of assertion (c). Here $\g$ is a
reductive Lie algebra and $\h$ is its indecomposable essential not
semisimple subalgebra. The proof is carried out in several steps.

{\it Step 1.} It follows from Proposition~\ref{Prop:4.9} that
$(\g,\h^{sat})$ is the direct sum of the pairs indicated in (c1).
(c2) is checked case by case.

{\it Step 2.} Let $\widetilde{\h}$ be an initial subalgebra of $\g$
and $\z=\z(\widetilde{\h})$.
 By the explicit
description of possible $\widetilde{\h}$ given above, we see that
$\dim\a(\g,[\widetilde{\h},\widetilde{\h}])=\dim\a(\g,\widetilde{\h})+\dim\z(\widetilde{\h})$.
It follows from the Frobenius reciprocity and the isomorphism
$V(\lambda)^{[\widetilde{\h},\widetilde{\h}]}\cong
V(\lambda^*)^{[\widetilde{\h},\widetilde{\h}]*}$ that
$\a(\g,[\widetilde{\h},\widetilde{\h}])$ is spanned by highest
weights $\lambda$ such that
$V(\lambda)^{[\widetilde{\h},\widetilde{\h}]}\neq 0$. Denote by
$\alpha_x(\lambda)$ the scalar by which $x\in
\z(\widetilde{\h})\cong
\widetilde{\h}/[\widetilde{\h},\widetilde{\h}]$ acts on
$K[G/(\widetilde{H},\widetilde{H})]^U_{(\lambda^*)}$.  By the
Frobenius reciprocity, this definition coincides with that given in
Remark~\ref{Rem:1.5}. It follows directly from definition that
$\alpha_x$ is a linear function on
$\a(\g,[\widetilde{\h},\widetilde{\h}])$ annihilating
$\a(\g,\widetilde{\h})$ and depending linearly on $x$.

{\it Step 3.} Now let $\widehat{\z}\subset \z:=\z(\widetilde{\h})$
be some algebraic subalgebra and
$\h_0=[\widetilde{\h},\widetilde{\h}]\oplus \widehat{\z}$. The pair
$(\g,\h_0)$ is indecomposable iff $\widehat{\z}$ projects nonzero to
any ideal of $\g$.  The characterization of essential subalgebras
with a fixed saturation given in the end of
Section~\ref{SECTION_essential} implies that $\h_0\subset \g$ is
essential if $(\g,\h_0)$ is indecomposable. This proves (c3).

{\it Step 4.} Let us check (c4). Note that $\h^{sat}=[\h,\h]\oplus
\z$, where $\z=\z(\z_\g([\h,\h]))$. The map $x\mapsto \alpha_x:
\z(\z_\g([\h,\h]))\rightarrow (\a(\g,[\h,\h])/\a(\g,[\h,\h]\oplus
\z))^*$ constructed on step 2 is injective by
Proposition~\ref{Prop:4.7}. Since $\dim \z=\dim
\a(\g,[\h,\h])/\a(\g,[\h,\h]\oplus\z)$ (see step 2), we see that
$x\mapsto \alpha_x$ is an isomorphism.

It remains to compute the functions $\alpha_x\in \a_{\g,[\h,\h]}^*$.
Let $\g=\z(\g)\oplus\g_1\oplus\ldots\oplus\g_k$ be the decomposition
of $\g$ into the direct sum of the center and the simple
non-commutative ideals,
$\h_i=[\h,\h]\cap\g_i,\widetilde{\h}_i=\h^{sat}\cap\g_i$. By step 1,
$[\h,\h]=\bigoplus_i\h_i, \h^{sat}=\bigoplus_i\widetilde{\h}_i$,
$\dim\widetilde{\h}_i/\h_i=1$. Clearly, to compute $\alpha_x$ it is
enough to assume $\g=\g_i,\h=\h_i$ (notice that $\h_i$ is not
necessarily essential subalgebra of $\g_i$ but it is does not
matter). In all cases  $\h$ is contained in the annihilator of the
highest weight of $V(\lambda_i)$, where $\lambda_i$ is contained in
column 3 in Table~\ref{Tbl:1.6}. Now the content of column 4 is
obtained by the direct computation. The spaces
$\a(\g_i,\widetilde{\h}_i) $ are computed in~\cite{Kramer} in all
cases except $(\g_i,\h_i)=(\sl_{n},\sl_k)$. To compute
$\a(\g_i,\widetilde{\h}_i)$ in this case it is enough to note that
$\pi_{n-k}^\vee\in\widetilde{\h}_i$  multiplies  an appropriate
highest vector in the $\g_i$-module $V(\pi_i)$ by $\frac{ik}{n}$ for
$i\leqslant n-k$ and by $\frac{(i-n)k}{n}$ for $i\geqslant k$.

{\it Step 5.} It remains to show (c5). Let $\widetilde{H},H$ be the
connected subgroups of $G$ corresponding to the algebras
$\h^{sat},\h$. These groups satisfy condition (1) of
Proposition~\ref{Prop:4.7}. It can be seen directly from the
definition of $\alpha_x$ that $x\in \widetilde{\h}/\h$ acts on
$V(\lambda)^{\h}$ by the multiplication by
$\alpha_{\widetilde{x}}(\lambda)$, where $\widetilde{x}$ is an
arbitrary element of $\z$ projecting to $x$. In particular,
$\alpha_x$ annihilates $\a(\g,\h)$ for $x\in \z(\h)$. Thus
$\a(\g,\h)/\a(\g,\h^{sat})$ lies in the annihilator of $\z(\h)$ in
$\a(\g,[\h,\h])/\a(\g,\h^{sat})$. To complete the proof it remains
to note that
\begin{align*}
\dim \z(\h)+\dim \a(\g,\h)/\a(\g,\h^{sat})=\dim\z(\h)+\dim
\widetilde{\h}/\h=\dim\a(\g,[\h,\h])/\a(\g,\h^{sat}).
\end{align*}

This completes the proof of Theorem \ref{Thm:1.3}.

In the end of the section we give an application of the
classification of essential subalgebras.  There are many papers
dealing with the classification of reductive subalgebras
$\h\subset\g$ with small complexity (the complexity of a subalgebra
is, by definition, the complexity of the corresponding homogeneous
space).

The classification of spherical subalgebras in simple Lie algebras
was given in~\cite{Kramer}. The partial classification for not
simple algebras $\g$ was given in~\cite{Brion},\cite{Mikityuk}. The
classification was completed in~\cite{Yakimova}. The classification
of subalgebras with complexity 1 in simple Lie algebras was carried
out in~\cite{Panyushev4}. The general case is treated in~\cite{ACh}.

Clearly,  $c_G(G/H)=c_{(G,G)}((G,G)/(G,G)\cap H)$. Therefore, to
classify all reductive subalgebras with given complexity it is
enough to consider the case of semisimple $\g$.

\begin{Prop}\label{Prop:5.1}
Let $H$ be a reductive subgroup in semisimple group $G$. Then
$c_G(G/H)=(\dim\g+\dim\l_0-\rank_G(G/H))/2-\dim\h$, where $\l_0$ is
the s.s.g.p. for the $\h^{ess}$-module $\g/\h^{ess}$
\end{Prop}
\begin{proof}
This follows from (\ref{eq:compl_formula}).
\end{proof}

So the  classification of subalgebras with given complexity follows
from the classification of subalgebras with given dimension  and
essential part.

\begin{Rem}
In fact, in the classification of essential subalgebras we used only
Kr\"{a}mer's classification of spherical reductive subalgebras in
simple algebras.
\end{Rem}

\section{Conventions and notations}\label{SECTION_Not_Conv}
\begin{longtable}{p{5.5cm} p{10cm}}
$\a_{G,X}$& the Cartan subspace of an irreducible $G$-variety $X$
\\ $\Aut(\g)$ (resp., $\Int(\g)$)& the group of all (resp., inner)
automorphisms of a Lie algebra $\g$.
\\ $c_G(X)$& the complexity of a $G$-variety $X$.
\\ $(G,G)$ (resp., $[\g,\g]$)& the commutant of an algebraic group
(resp., a Lie algebra $\g$)
\\ $G^{\circ}$& the connected component of unit of an algebraic
group $G$.
\\ $G_x$& the stabilizer of a point $x\in X$ under an action
$G:X$.
%\\ $U^0$& the annihilator of a subspace $U\subset V$ in the dual
%space $V^*$.
\\ $N_G(H)$, (resp. $N_G(\h),\n_\g(\h)$)& the normalizer of an
algebraic subgroup $H$ in an algebraic group $G$ (resp. of a
subalgebra $\h\subset \g$  in an algebraic group $G$, of a
subalgebra $\h\subset \g$ in a Lie algebra $\g$).
\\ $\rank_G(X)$& the
rank of an irreducible $G$-variety $X$
\\ s.g.p. & stabilizer in general position.
\\ s.s.g.p.& stable subalgebra in general position.
\\ $V^\g$& $=\{v\in V| \g v=0\}$, where $\g$ is a Lie algebra and
$V$ is a $\g$-module.
\\ $V_{(\lambda)}$& the isotypical component of a $\g$-module
$V$ ($\g$ is a reductive Lie algebra) corresponding to a highest
weight $\lambda$.
\\ $V(\mu)$& the irreducible module with the highest weight $\mu$
over a reductive algebraic group  or a reductive Lie algebra.\\
$\X(G)$& the weight lattice of an algebraic group $G$\\  $\X_{G,X}$&
the weight lattice of an irreducible $G$-variety $X$\\    $X^G$& the
fixed-point set
for an action of $G$ on $X$\\
$\#X$& the number of elements in a finite set $X$.
\\ $\z(\g)$& the center of a Lie algebra $\g$.
\\   $Z_G(H)$, (resp. $Z_G(\h),\z_\g(\h)$)& the centralizer of an
algebraic subgroup $H$ in an algebraic group $G$ (resp. of a
subalgebra $\h\subset \g$  in an algebraic group $G$, of a
subalgebra $\h\subset \g$ in a Lie algebra $\g$).
\\ $\alpha^\vee$& the dual root corresponding to a root $\alpha$.
\\  $\Delta(\g)$& the root system of a reductive Lie algebra $\g$.
\\ $\lambda^*$& the dual highest weight to a highest weight
$\lambda$.
\end{longtable}

If an algebraic group is denoted by a capital Latin letter, then its
Lie algebra is  denoted by the corresponding small fraktur letter.

All homomorphisms   of reductive algebraic Lie algebras (for
example, representations) are assumed to be the differentials of
homomoprhisms of the corresponding reductive algebraic groups.

When $\g$ is a simple Lie algebra, then $\alpha_i$ denotes its
$i$-th simple root and $\pi_i$ denotes the corresponding fundamental
weight. For roots and weights of simple Lie algebras we use the
notation taken from~\cite{VO}.

Now we explain notation for subalgebras in semisimple Lie algebras.
For subalgebras of exceptional Lie algebras we use the notation from
\cite{Dynkin}.

If $\g=\sl_n$, then by $\sl_k,\so_k,\sp_{k}$ we mean subalgebras,
that annihilate some subspace  of dimension $n-k$ in $\C^n$ and
preserve a corresponding form on its complement (for $\so_k,\sp_k$).
The subalgebras $\so_k\subset\so_n, \sp_k\subset \sp_n$ are defined
similarly. A subalgebra $\gl_k$ is embedded into $\so_n$ by the
direct sum of $\tau,\tau^*$ and the trivial representation (where
$\tau$ denotes the tautological representation of $\gl_k$). The
embeddings of $\sl_k,\sp_k$ into $\so_n$ are the compositions of the
described embedding $\gl_k\hookrightarrow \so_n$ and the embeddings
$\sp_k\hookrightarrow \gl_k, \sl_k\hookrightarrow \gl_k$. The
subalgebra $G_2$ (resp., $\spin_7$) of $\so_n$ is the image of
$G_2$,(resp., $\so_7$) under the direct sum of the 7-dimensional
irreducible (resp., spinor)  and trivial representations.

The description above determines the subalgebras up to conjugacy in
$\Aut(\g)$.

\bigskip

{\Small Department of Higher Algebra, Department of Mechanics and
Mathematics, Moscow State University.

E-mail address: ivanlosev@yandex.ru}
\end{document}